\documentclass[a4paper,11pt]{article}

\usepackage{ucs}
\usepackage[utf8]{inputenc}

\usepackage{graphicx}
\usepackage{amsfonts}
\usepackage{dsfont}
\usepackage{amssymb}
\usepackage{amsmath}
\usepackage{amsthm}
\usepackage{enumerate}
\usepackage{stmaryrd}
\usepackage{fullpage}
\usepackage{ifthen}
\usepackage{subfigure}
\usepackage{epic}
\usepackage{authblk}
\usepackage{textcomp}
\usepackage{mathrsfs}
\usepackage{bm}
\usepackage[small]{caption}
\usepackage{tikz,tkz-tab}
\usetikzlibrary{matrix}
\usepackage{pgfplots}
\pgfplotsset{compat=newest} 
\pgfplotsset{plot coordinates/math parser=false}

\usepackage[hypertexnames=false,colorlinks=true,linkcolor=blue,citecolor=blue]{hyperref}
\usepackage[numbers,comma,square,sort&compress]{natbib}
\usepackage[letterpaper,text={6.5in,9in},centering]{geometry}

\usepackage{color}
\usepackage{titlesec}

\graphicspath{{eps/}{pdf/}{png/}}

\newcommand{\bqq}{\begin{equation}}
\newcommand{\eqq}{\end{equation}}
\newcommand{\bqs}{\begin{equation*}}
\newcommand{\eqs}{\end{equation*}}

\newcommand{\R}{\mathbb{R}} 
\newcommand{\N}{\mathbb{N}}

\newcommand{\Z}{\mathbb{Z}}

\newcommand{\E}{\mathcal{E}}
\newcommand{\V}{\mathcal{V}}

\newcommand{\G}{\mathcal{G}}

\newcommand{\md}{\mathrm{d}}

\newtheorem{lem}{Lemma}[section]
\newtheorem{thm}{Theorem}
\newtheorem{prop}[lem]{Proposition}

 {\begin{trivlist}\item[]\textbf{Proof#1 }}%
 {\hspace*{\fill}$\rule{0.3\baselineskip}{0.35\baselineskip}$\end{trivlist}}

  {\begin{trivlist}\item[]{\bf Hypothesis #1 }\em}{\end{trivlist}}

\numberwithin{equation}{section}

\title{Spreading properties for SIR models on homogeneous trees}

\author[1]{Christophe Besse \& Gr\'egory Faye\footnote{Corresponding author: \texttt{gregory.faye@math.univ-toulouse.fr}}}
\affil[1]{\small CNRS, UMR 5219, Institut de Math\'ematiques de Toulouse, 31062 Toulouse Cedex, France}

\begin{document}
\maketitle

\begin{abstract}
We consider an epidemic model of SIR type set on a homogeneous tree and investigate the spreading properties of the epidemic as a function of the degree of the tree, the intrinsic basic reproduction number and the strength of the interactions within the population of infected individuals. When the degree is one, the homogeneous tree is nothing but the standard lattice on the integers and our model reduces to a SIR model with discrete diffusion for which the spreading properties are very similar to the continuous case. On the other hand, when the degree is larger than two, we observe some new features in the spreading properties. Most notably, there exists a critical value of the strength of interactions above which spreading of the epidemic in the tree is no longer possible.
\end{abstract}

\bigskip

\section{Introduction}

Let $\G=(\V,\E)$ be a connected graph where $\V$ is the set of vertices and $\E$ is the set of edges. We introduce the following SIR model (Kermack and McKendrick \cite{KMK27}) on the graph $\G$
\bqq
\left\{
\begin{split}
S_v'(t)&= -\tau S_v(t)I_v(t),\\
I_v'(t)&= \tau S_v(t)I_v(t)-\eta I_v(t)+ \lambda  \sum_{v'\sim v} \left( I_{v'}(t)-I_v(t)\right),\\
R_v'(t)&=\eta I_v(t),
\end{split}
\right. \quad v\in\V, \quad t>0,
\label{SIRgraph}
\eqq
where $S_v$ stands for the density of susceptible individuals, $I_v$ represents the density of infected individuals and $R_v$ is the density of removed individuals at vertex $v$. Here, all parameters $\tau>0$, $\eta>0$ and $\lambda>0$ are set to be positive and homogeneous, in the sense that they do not depend on the vertex $v$. In this setting, $\tau>0$ is a contact rate between susceptible and infected populations while $1/\eta>0$ is the average infectious period. We refer to  \cite{Heth00} for a review on SIR models.

Exchanges of infected individuals in the graph are modeled by the term $\lambda \sum_{v'\sim v} \left( I_{v'}(t)-I_v(t)\right)$ with $\lambda>0$ the strength of the exchanges and where the sum is taken on all adjacent vertices $v'$ to vertex $v$ which we denote with the shorthand notation $v'\sim v$. By adjacent, we mean that there exists an edge $e\in\E$ such that $e$ connects $v$ to $v'$. In graph theory \cite{bollobas}, the term $\sum_{v'\sim v} \left( I_{v'}(t)-I_v(t)\right)$ is often referred to as the graph Laplacian of $\G$ evaluated at vertex $v$ and models diffusion of infected individuals within the graph and accordingly the parameter $\lambda$ can be interpreted as a diffusion coefficient. SIR models on graphs with similar exchanges of infected individuals have previously been introduced in the literature \cite{BH13,CHS18,FGW16,CFG06,CGH17} with a particular emphasis on spreading properties. Actually, model \eqref{SIRgraph} is at the crossroad of traditional reaction-diffusion models on graphs \cite{BCVV12,HH19} and epidemic models that incorporate more sophisticated interactions dynamics \cite{BDL08,SSVM13,HSMS13,SB19,SCB15,BB20,BG18,BF21}. Note that in \eqref{SIRgraph}, we have assumed that only the infected population is subject to diffusion within the graph, and we think of $S_v$ being an ambient population whose movement does not affect its distribution. This is of course a strong biological limitation and considering the susceptible population as an ambient population is a first step. It would be natural to extend our model to the case that individuals in the susceptible population can also interact through exchanges of the form $\sum_{v'\sim v} \left( S_{v'}(t)-S_v(t)\right)$. In the case of continuous spatially extended systems of reaction-diffusion type, it is notorious that allowing the susceptible population to diffuse is more challenging from a theoretical point of view as monotonicity properties of the solutions are lost \cite{BNR20} and we leave such an analysis for a future work. 

In this work, we will focus on the spreading properties of system \eqref{SIRgraph}. Namely, we would like to characterize the long time dynamics of the solutions of \eqref{SIRgraph} starting from the initial configuration where susceptible individuals are homogeneously distributed across the graph, that is $S_v(t=0)=s_0\in(0,1)$ for each $v\in\V$, and where infected populations are only present at finitely many vertices. To simplify, we will sometimes consider the case where infected individuals are initially present at only one given vertex and without loss of generality we shall always assume that $R_v(t=0)=0$ for all $v\in\V$. In what follows, we will ignore the dynamics on $R_v(t)$ since it can be read out from $I_v(t)$ via $R_v(t)=\eta \int_0^t I_v(s)\md s$ for each $v\in\V$ and $t>0$. Assuming a homogeneous distribution across the graph of susceptible individuals is questionnable from a biological point of view as in practical situations this distribution is most likely to be heterogeneous. Here, we adopt this formalism since it will allow us to carry a fairly complete mathematical analysis with closed form formulas which are relatively simple to interpret. Moreover, the homogeneous case already sheds light on the effects of networks structure on the propagation of epidemics.

\begin{figure}[t!]
  \centering
  \includegraphics[width=.5\textwidth]{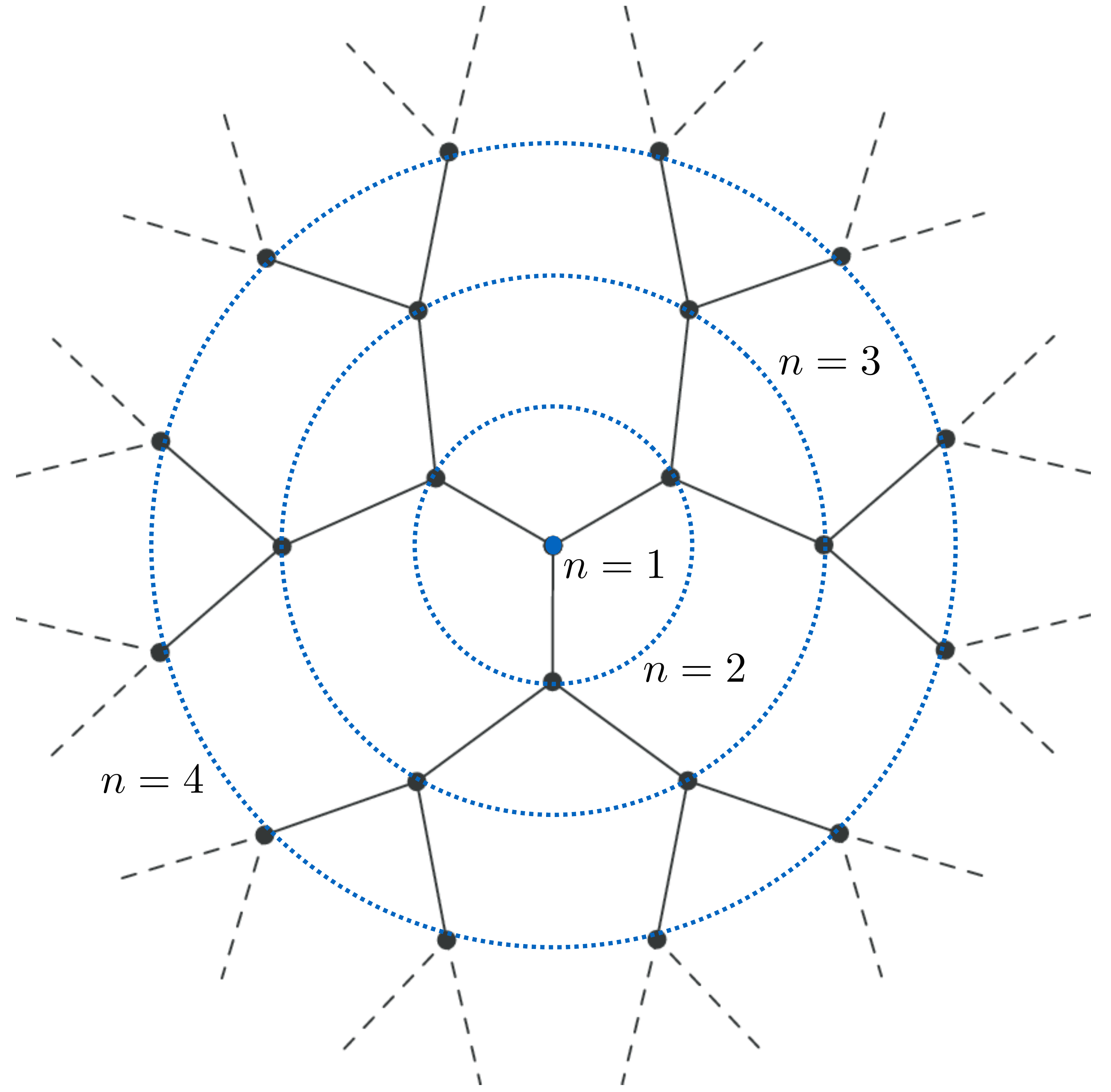} 
  \caption{Example of the homogeneous tree $\mathbb{T}_2$ of degree $k=2$ where each node has precisely $3$ adjacent vertices. Here, we have identified one vertex has being the root of the tree and denoted it $n=1$. In model \eqref{SIRtree}, we label $(S_n,I_n)$ as a representative vertex from the set of vertices at distance $n-1$ from the root. For example at distance $2$ from the root there are $6$ vertices which are all identified. Infected individuals are initially present only at the root while susceptible individuals are homogeneously distributed across the tree.}
  \label{fig:tree}
\end{figure}

Coming back to the long time dynamics of the solutions to system \eqref{SIRgraph}, we would like to determine under which conditions on the parameters and on the initial configuration an epidemic may spread in the graph, and characterize at which speed this spreading occurs and what will be the final configuration. Without further assumption on the graph $\G$, it is very difficult to provide any answer to the above questions. This is why in this work, we will focus on the specific case where the graph $\G$ is a homogenous tree of degree $k \in \N$ with $k\geq1$, which we will denote $\mathbb{T}_k$ from now on\footnote{Note that in the physics literature homogeneous trees of degree $k$ are often called Bethe lattices \cite{bethe}.}. A homogeneous tree of degree $k$ is an infinite graph where each vertex $v$ has precisely $k+1$ adjacent vertices. We refer to Figure~\ref{fig:tree} for an illustration in the case $k=2$. 

When $k=1$, the homogeneous tree $\mathbb{T}_1$ is nothing but $\Z$ the lattice of the integers. In that case, system \eqref{SIRgraph} reduces to
\bqq
\left\{
\begin{split}
S'_j(t)&=-\tau S_j(t)I_j(t),\\
I'_j(t)&=\tau S_j(t)I_j(t)-\eta I_j(t)+\lambda\left( I_{j-1}(t)-2I_j(t)+I_{j+1}(t)\right),
\end{split}
\right. 
\label{SIRlattice}
\eqq
for $j\in\Z$ and $t>0$. The graph Laplacian at vertex $j$ given by the term $I_{j-1}(t)-2I_j(t)+I_{j+1}(t)$ takes the traditional form of a discrete Laplacian. Actually, system \eqref{SIRlattice} can be interpreted as a discretized version, through finite differences, of the following spatially continuous SIR model set on $x\in\R$
\bqq
\left\{
\begin{split}
\partial_t S(t,x)&=-\tau S_j(t,x)I(t,x),\\
\partial_t I(t,x)&=\tau S_j(t,x)I(t,x)-\eta I(t,x)+ d \partial_x^2 I(t,x),
\end{split}
\right. 
\label{SIRcontinuous}
\eqq
where $d>0$ is some diffusion coefficient by setting $\lambda = \frac{d}{\Delta x^2}$ for some small $\Delta x$. The continuous model \eqref{SIRcontinuous} has received much attention in the past decades and especially its spreading properties, see for example  \cite{A77,BRR21,BNR20,DG14,W82} and references therein. One of our objective will be to understand how these spreading properties, which hold true in the continuous case, will persist (or not) in our discrete setting. It turns out, that the degree $k$ of the homogeneous tree $\mathbb{T}_k$ will be a key parameter and phenomenologically new behaviors will emerge for large values of $k$. 

When $k\geq2$, we will consider initial condition where the density of infected individuals is non zero at only one vertex. By convention, we will label this vertex as $1$ with associated density $(S_1,I_1)$, and it will be identified as the root of the tree. By symmetry of the model, all densities $(S_v,I_v)$ at some fixed distance away from the root are equal. As a consequence, it will be convenient to denote $(S_n,I_n)$ as a representative vertex from the set of vertices at distance $n-1$ from the root, see Figure~\ref{fig:tree} for an illustration in the case $k=2$. With these notations, system \eqref{SIRgraph} becomes
\bqq
\left\{
\begin{split}
S'_n(t)&=-\tau S_n(t)I_n(t),\\
I'_n(t)&=\tau S_n(t)I_n(t)-\eta I_n(t)+\lambda\left(I_{n-1}(t)-(k+1)I_n(t)+kI_{n+1}(t)\right),
\end{split}
\right. 
\label{SIRtree}
\eqq
for $n\geq2$ and
\bqq
\left\{
\begin{split}
S'_1(t)&=-\tau S_1(t)I_1(t),\\
I'_1(t)&=\tau S_1(t)I_1(t)-\eta I_1(t)+\lambda(k+1)\left(-I_1(t)+I_2(t)\right),
\end{split}
\right. 
\label{SIRroot}
\eqq
for the equation at the root. We remark that for $n\geq2$, the diffusive term can be expressed as 
\bqs
I_{n-1}(t)-(k+1)I_n(t)+kI_{n+1}(t)=I_{n-1}(t)-2I_n(t)+I_{n+1}(t) +(k-1)(I_{n+1}(t)-I_n(t)).
\eqs
As previously noticed in a different context  \cite{HH19}, system \eqref{SIRtree}-\eqref{SIRroot} can also be interpreted as a discretization of the following continuous model for $x>0$
\bqq
\left\{
\begin{split}
\partial_t S(t,x)&=-\tau S_j(t,x)I(t,x),\\
\partial_t I(t,x)&=\tau S_j(t,x)I(t,x)-\eta I(t,x)+d \partial_x^2 I(t,x)+\frac{1}{\epsilon}(k-1)\partial_x I(t,x),
\end{split}
\right. 
\label{SIRcontinuousMod}
\eqq
and no-flux boundary condition at the left boundary $x=0$ by setting once again $\lambda = \frac{d}{\Delta x^2}$  and $\epsilon = \frac{\sqrt{\Delta x}}{d}>0$ for some small $\Delta x$. The main difference is the new drift term which is penalized by $\frac{1}{\epsilon}$. One expects that there will be a trade-off where either this advection dominates the dynamics and compactly supported initial conditions for \eqref{SIRcontinuousMod} are propagated to the left of the domain and eventually converge to zero, or reaction terms dominate and compactly supported initial conditions will spread across the domain. One of our objectives is precisely to understand the possible transition from an epidemic spreading
to pointwise convergence to zero and, in the case of spreading, to predict the spreading speed of the solution as a function of the parameters \eqref{SIRtree}-\eqref{SIRroot}. Let us finally note that similar behavior have been described for reaction-diffusion equation of Fisher-KPP type set on hyperbolic spaces \cite{MPT15}.

\section{Main results}

In this section, we present our main results and distinguish between the case of the lattice $\Z$ and a homogeneous tree $\mathbb{T}_k$ of degree $k\geq2$.

\subsection{Case of the lattice $\Z$}
\begin{figure}[t!]
  \centering
  \includegraphics[width=.5\textwidth]{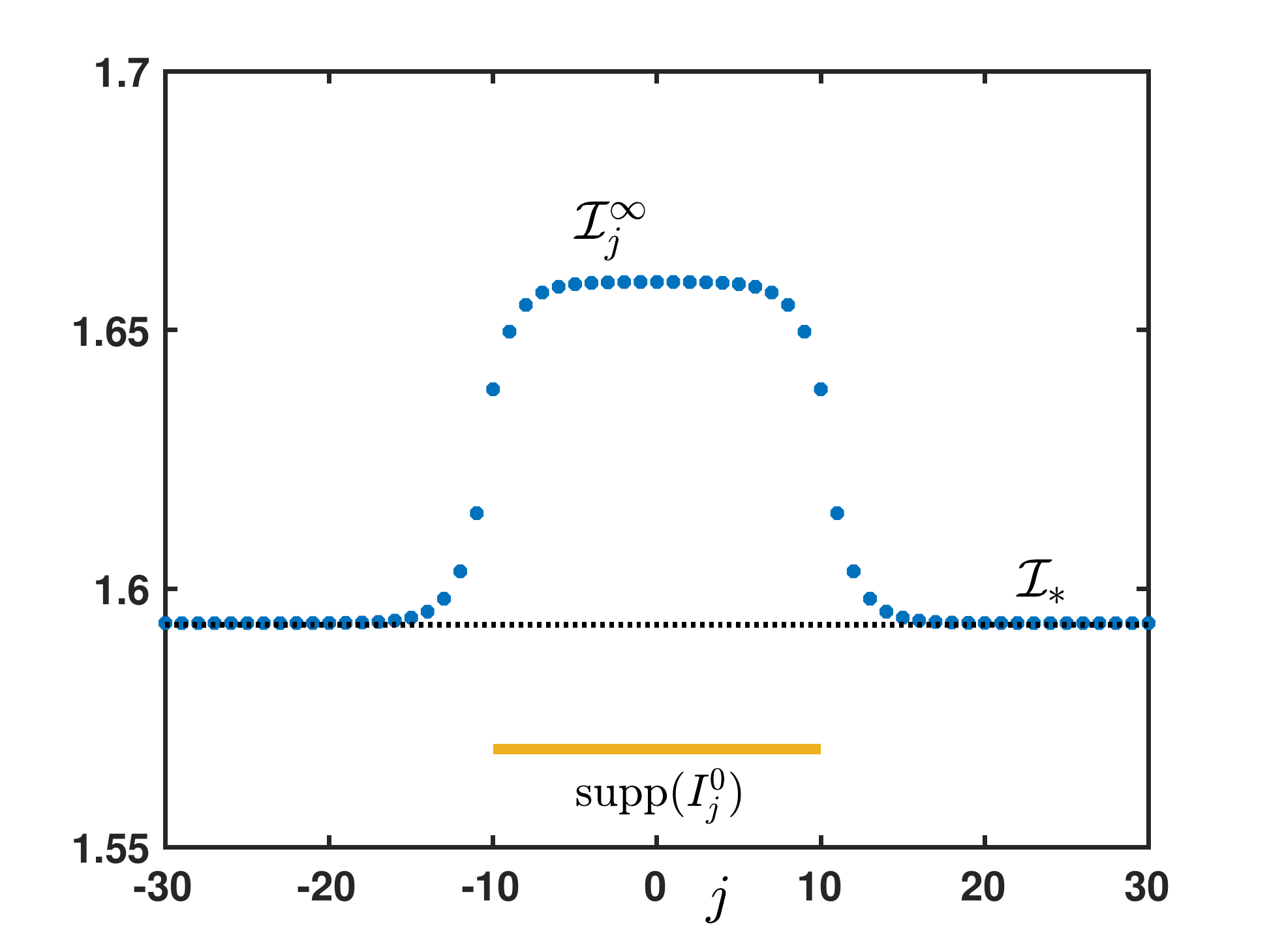} 
  \caption{Numerically computed stationary solution $\left(\mathcal{I}_j^\infty\right)_{j\in\Z}$ in the case $\mathscr{R}_0>1$ for a given $\lambda>0$. The support of the initial density of infected individuals is represented by the yellow bar and was set to $\llbracket-10,10 \rrbracket$.}
  \label{fig:Iinf}
\end{figure}

We complement \eqref{SIRlattice} with an initial condition of the form
\bqq
S_j(t=0)=s_0, \quad I_j(t=0)=I_j^0 , \quad j\in\Z,
\label{IC}
\eqq
where $s_0\in(0,1)$ and $I_j^0\in(0,1)$ for each $j\in\Z$ has finite support.  As already stated in the introduction, our aim is to investigate the long time dynamics of \eqref{SIRlattice} subject to the initial condition \eqref{IC}. We define the cumulative density of infected individuals at lattice site $j$ and time $t$ as $\mathcal{I}_j(t):=\int_0^tI_j(s)\md s$ such that the density of susceptible individuals can be expressed as
\bqs
\ln\left(\frac{S_j(t)}{s_0}\right)=-\tau \mathcal{I}_j(t).
\eqs
Thus $\mathcal{I}_j(t)$ satisfies the lattice differential equation
\bqq
\mathcal{I}_j'(t)=f(\mathcal{I}_j(t))+I_j^0+\lambda\left( \mathcal{I}_{j-1}(t)-2\mathcal{I}_j(t)+\mathcal{I}_{j+1}(t)\right),
\label{KPPlike}
\eqq
together with the initial condition 
\bqq
\mathcal{I}_j(t=0)=0, \quad\text{ for all }\quad j\in\Z.
\label{ICkpp}
\eqq
Here, the nonlinearity $f$ is given by
\bqs
f(v):=s_0\left(1-e^{-\tau v}\right)-\eta v,
\eqs
which is smooth, concave on $[0,+\infty)$ and vanishes at $v=0$. For future reference, we note that
\bqs
f'(0)=\eta \left( \mathscr{R}_0-1\right), \quad \mathscr{R}_0:= \frac{s_0\tau}{\eta}.
\eqs
The quantity $\mathscr{R}_0$ is the usual basic reproduction number \cite{DHM90,VanW02}. Our first result states the existence of a unique positive, bounded, stationary solution to \eqref{KPPlike} and characterizes its asymptotic behavior, see Figure~\ref{fig:Iinf} for an illustration.

\begin{thm}\label{thm1}
The equation \eqref{KPPlike} admits a unique positive, bounded, stationary solution $\left(\mathcal{I}_j^\infty\right)_{j\in\Z}$ which satisfies
\bqs
\underset{|j|\rightarrow+\infty}{\lim}\mathcal{I}_j^\infty=\left\{
\begin{array}{lcl}
0, & \text{if}  & \mathscr{R}_0 \leq 1, \\ 
\mathcal{I}_*, & \text{if} & \mathscr{R}_0>1,
\end{array}
\right.
\eqs
where $\mathcal{I}_*>0$ is the unique positive zero of $f$.
\end{thm}

It is interesting to point out that \eqref{KPPlike} can be interpreted as a reaction-diffusion equation of Fisher-KPP type on the lattice $\Z$ with a heterogeneity given by $I_j^0$. We refer to \cite{BCVV12,HH19} for spreading properties of the Fisher-KPP equation set on the lattice $\Z$. In our SIR epidemic setting, as in the continuous case \cite{BRR21}, it turns out that the above stationary solution is a global attractor for the dynamics of \eqref{KPPlike} which is the result of the following theorem.

\begin{thm}\label{thm2}
Let $\left(\mathcal{I}_j(t)\right)_{j\in\Z}$ be the solution of \eqref{KPPlike} starting from some nonnegative bounded compactly supported initial condition.
Then  $\left(\mathcal{I}_j(t)\right)_{j\in\Z}$ converges as $t\rightarrow+\infty$ locally uniformly to $\left(\mathcal{I}_j^\infty\right)_{j\in\Z}$.
\end{thm}

A consequence of the above Theorem~\ref{thm2} is that the convergence also holds true for the time derivative of $\left(\mathcal{I}_j(t)\right)_{j\in\Z}$ such that we can deduce that locally uniformly in $j$ we have
\bqs
I_j(t)=\mathcal{I}_j'(t) \rightarrow 0 \text{ as } t\rightarrow +\infty,
\eqs
which means that the density of infected individuals asymptotically vanishes in time at each fixed lattice site. Recalling that the density of susceptible individuals can be read out from $\mathcal{I}_j(t)$ via
\bqs
S_j(t)=s_0e^{-\tau \mathcal{I}_j(t)},
\eqs
we can also quantify the density of individuals that will be infected during the course of the epidemic at a given lattice site $j$, that we denote $\mathcal{I}_j^{tot}$. It is given by
\bqs
\mathcal{I}_j^{tot}=s_0\left(1-e^{-\tau \mathcal{I}_j^\infty}\right), \quad j\in\Z.
\eqs
As in the fully continuous setting \cite{DG14,BRR21}, we remark that $\left(\mathcal{I}_j^{tot}\right)_{j\in\Z}$ is not constant, which comes from our assumption that susceptible individuals do not diffuse on the lattice. Now, using the result of Theorem~\ref{thm1}, we obtain the following dichotomy 
\bqs
\underset{|j|\rightarrow+\infty}{\lim}\mathcal{I}_j^{tot}=\left\{
\begin{array}{lcl}
0, & \text{if}  & \mathscr{R}_0 \leq 1, \\ 
s_0\left(1-e^{-\tau \mathcal{I}_*}\right), & \text{if} & \mathscr{R}_0>1.
\end{array}
\right.
\eqs
This basically says that when $\mathscr{R}_0 \leq 1$ the epidemic does not propagate across the network and lattice sites very far from the support of the initial condition will effectively be not infected. On the other hand, when the basic reproduction number verifies $\mathscr{R}_0>1$, the epidemic propagates everywhere across the lattice. Very far from the support of the initial condition, there will be a portion $1-e^{-\tau \mathcal{I}_*}$ of infected individuals of the overall population. When $\mathscr{R}_0>1$, we can further characterize at which speed the epidemic spreads into the lattice and this is precisely the result of the next theorem. 

\begin{thm}\label{thm3}
Assume that $\mathscr{R}_0>1$. We define $c_*>0$ as
\bqs
c_*:=\underset{\gamma>0}{\min}~  \frac{\eta \left( \mathscr{R}_0-1\right)+\lambda \left( e^{-\gamma}-2+e^{\gamma} \right)}{\gamma}.
\eqs
Then, the solution $\left(\mathcal{I}_j(t)\right)_{j\in\Z}$ of \eqref{KPPlike}-\eqref{ICkpp} satisfies:
\begin{itemize}
\item[(i)] $\forall c\in(0,c_*)$,
\bqs
\underset{t\rightarrow+\infty}{\limsup} \left( \underset{|j| \leq ct}{\sup}\left| \mathcal{I}_j(t)-\mathcal{I}_j^\infty \right| \right) =0;
\eqs
\item[(ii)] $\forall c>c_*$,
\bqs
\underset{t\rightarrow+\infty}{\limsup} \left( \underset{|j| \geq ct}{\sup}\left| \mathcal{I}_j(t) \right| \right) =0.
\eqs
\end{itemize}
\end{thm}
The quantity $c_*$ is the asymptotic speed of spreading of the epidemic wave and it coincides with the asymptotic speed of spreading for the Fisher-KPP equation set on the lattice \cite{HH19}. As already noticed in \cite{HH19}, we have the following asymptotic of the spreading speed $c_*$ when $\lambda$ is small
\bqs
c_* \sim W_0\left( \frac{\eta(\mathscr{R}_0-1)}{\lambda}\right) \text{ as } \lambda \rightarrow 0,
\eqs
where $W_0$ denotes the principal branch of the Lambert W function, that is the multivalued inverse relation of the function $f(w) = we^w$ for $w\in 
\mathbb{C}$ \cite{lambert}. On the other hand, when $\lambda$ is large, we have that
\bqs
c_* \sim 2\sqrt{\eta(\mathscr{R}_0-1)\lambda} \text{ as } \lambda \rightarrow +\infty.
\eqs
Let us remark that $2\sqrt{\eta(\mathscr{R}_0-1)\lambda}$ is precisely the asymptotic speed of propagation in the continuous case \cite{BRR21}. It turns out that the speed $c_*$ also characterizes the threshold for the existence of traveling wave solutions associated to system \eqref{SIRlattice}. For a given $s_0\in(0,1)$, by traveling wave solution to system \eqref{SIRlattice} we mean a solution of \eqref{SIRlattice} which takes the form $(S_j(t),I_j(t))=(S(j-ct),I(j-ct))$ with profiles $(S(x),I(x))$ and wave speed $c$ solutions of
\bqq
\left\{
\begin{split}
-cS'(x)&=-\tau S(x)I(x),\\
-cI'(x)&=\tau S(x)I(x)-\eta I(x)+\lambda\left( I(x-1)-2I(x)+I(x+1)\right),\\
c>0, &\quad I>0 \text{ is bounded }, \quad 0<S<s_0,
\end{split}
\right.
\label{TWlattice}
\eqq
with asymptotic conditions 
\bqq
\left\{
\begin{split}
S(+\infty)&=s_0,\\
I(\pm\infty)&=0.
\end{split}
\right.
\label{TWasymp}
\eqq
Here the wave speed $c>0$ is to be determined and we do not impose a prescribed value for $S$ at $-\infty$.

\begin{thm}\label{thm4}
Assume that $\mathscr{R}_0>1$. Then there exist traveling wave solutions with speed $c$ for any $c\geq c_*$. In addition,
\bqs
S(-\infty)=s_\infty>0,
\eqs
where $s_\infty$ is the unique positive real such that $\Psi(s_\infty)=\Psi(s_0)$ with $0<s_\infty<s_0$ and $\Psi(v):=v-\frac{\eta}{\tau}\ln(v)$. We further have that $s_\infty=s_0e^{-\tau \mathcal{I}_*}$.
For $c\in(0,c_*)$, no such traveling wave solutions exist. 
\end{thm}

Existence and uniqueness results of traveling wave solutions for discrete epidemic models on a lattice have been obtained recently \cite{FGW16,Wu17,CGH17}, but not directly for system \eqref{TWlattice}. Indeed, in \cite{FGW16,Wu17,CGH17}, the density of susceptible individuals is also subject to discrete diffusion which prevents us to directly apply to our system the aforementioned  results. Nevertheless, we can easily notice that system \eqref{TWlattice}-\eqref{TWasymp} can be recast as a traveling wave problem for the discrete Fisher-KPP problem by introducing $\mathcal{I}(x):= \frac{1}{c}\int_x^{+\infty} I(z)\md z$ for which existence and uniqueness results are readily available \cite{ZHH91,CG02,CFG06}.

\begin{figure}[t!]
  \centering
  \includegraphics[width=.5\textwidth]{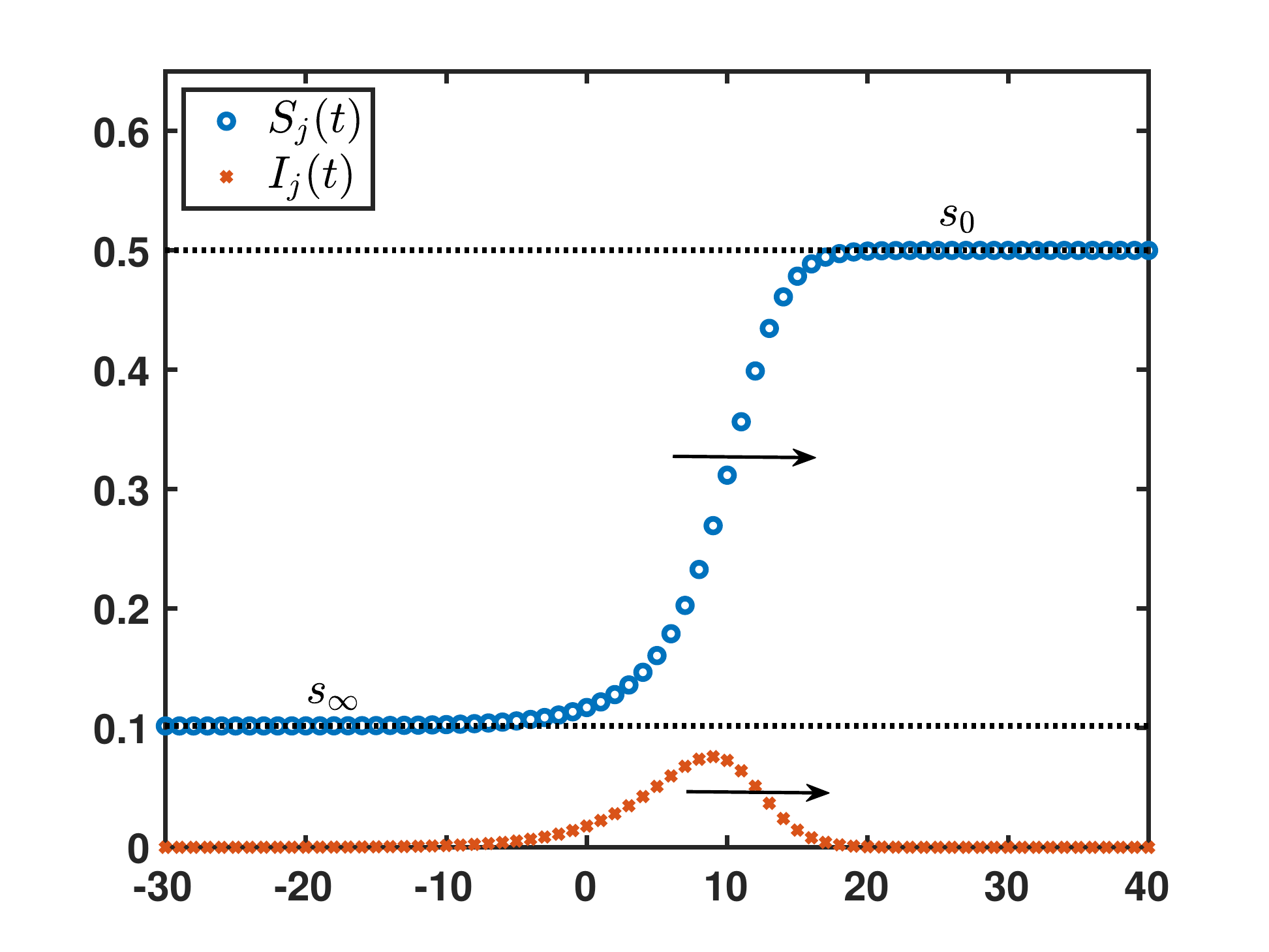} 
  \caption{Illustration of the traveling wave profiles $(S_j(t),I_j(t))=(S(j-ct),I(j-ct))$ from Theorem~\ref{thm4}.}
  \label{fig:TW}
\end{figure}

\subsection{Case of a homogeneous tree of degree $k\geq 2$}

\begin{figure}[t!]
\centering
\includegraphics[height=0.4\textheight]{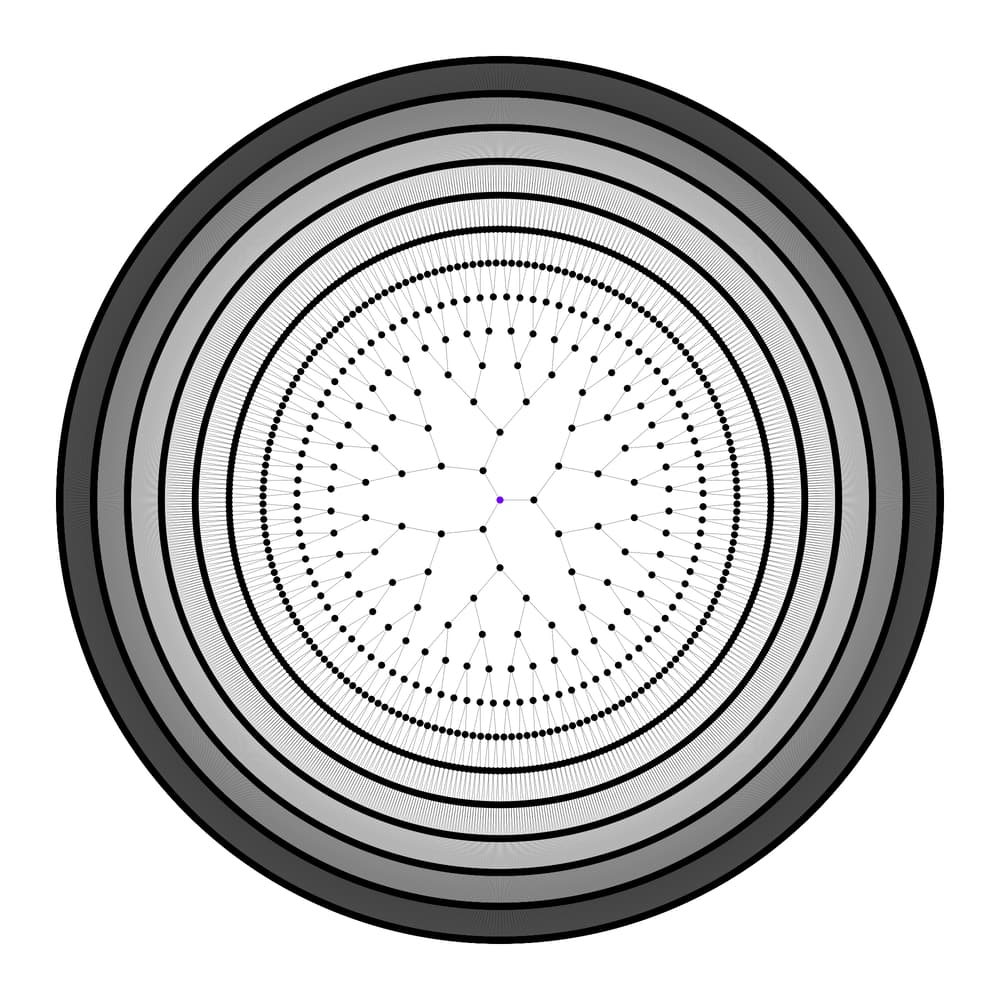}
\caption{Illustration of the initial condition for the density of infected individuals on the homogeneous tree of degree $k=2$. Infected individuals are only present at the root (purple dot in the center) with $I_1(t=0)=i_0\in(0,1)$ and $I_n(t=0)=0$ for $n\geq2$. Here, we represented the tree up to $n=14$ generations.}
\label{fig:InitialCond}
\end{figure}

We now turn our attention to the case where our epidemic model \eqref{SIRgraph} is set on a homogeneous tree of degree $k\geq 2$. Throughout this section, we will  assume that the initial density of infected individuals is non zero at only one vertex of the graph. We recall from the introduction that we label this vertex as $1$ and denote $(S_1(t),I_1(t))$ the associated density. Using the symmetry of the model, we also denote $(S_n(t),I_n(t))$ a representative vertex from the set of vertices at distance $n-1$ from the root. As a consequence, the sequence $(S_n(t),I_n(t))_{n\geq1}$ is solution of the system \eqref{SIRtree}-\eqref{SIRroot} which is complemented with an initial of the form
\bqs
S_n(t=0)=s_0, \quad n\geq 1, \quad I_1(t=0)=i_0, \quad I_n(t=0)=0, \quad n\geq2,
\eqs
for some $s_0\in(0,1)$ and $i_0\in(0,1)$. We refer to Figure~\ref{fig:InitialCond} for an illustration of the configuration of the initial condition for the density of infected individuals in the case $k=2$.

Once again, for each $n\geq1$, we define the cumulated density of infected individuals as $\mathcal{I}_n(t):=\int_0^t I_n(s)\md s$ which satisfies
\bqq
\left\{
\begin{split}
\mathcal{I}'_n(t)&=f( \mathcal{I}_n(t))+\lambda\left(\mathcal{I}_{n-1}(t)-(k+1)\mathcal{I}_n(t)+k\mathcal{I}_{n+1}(t)\right), \quad n \geq 2,\\
\mathcal{I}'_1(t)&=f(\mathcal{I}_1(t))+i_0+\lambda(k+1)\left(-\mathcal{I}_1(t)+\mathcal{I}_2(t)\right),
\end{split}
\right. 
\label{KPPtree}
\eqq
together with the initial condition
\bqq
\mathcal{I}_n(t=0)=0, \quad\text{ for all }\quad n\geq 1.
\label{KPPtreeIC}
\eqq

We start by studying stationary solutions to \eqref{KPPtree}, and as in Theorem~\ref{thm1} on the lattice, we have existence and uniqueness of a positive bounded stationary solution. Similarly, when the basic reproduction number satisfies $\mathscr{R}_0 \leq 1$, the stationary solution $\mathcal{I}_n^\infty$ converges to zero as $n$ goes to infinity.

\begin{thm}\label{thm5}
The equation \eqref{KPPtree} admits a unique positive, bounded, stationary solution $\left(\mathcal{I}_n^\infty\right)_{n\geq 1}$. Furthermore, when $\mathscr{R}_0 \leq 1$, we have
\bqs
\underset{n\rightarrow+\infty}{\lim}\mathcal{I}_n^\infty=0.
\eqs
\end{thm}

It turns out that the asymptotic behavior of the stationary solution when $\mathscr{R}_0>1$ is more intricate than in the case on the lattice. We have the following result.

\begin{thm}\label{thm6}
Assume that $\mathscr{R}_0 > 1$. Then the unique stationary solution $\left(\mathcal{I}_n^\infty\right)_{n\geq 1}$ to \eqref{KPPtree} satisfies
\bqs
\underset{n\rightarrow+\infty}{\lim}\mathcal{I}_n^\infty=\left\{
\begin{array}{lcl}
\mathcal{I}_*, & \text{if}  & 0<\lambda<\lambda_c, \\ 
0, & \text{if} & \lambda>\lambda_c,
\end{array}
\right.
\eqs
where $\mathcal{I}_*>0$ is the unique positive zero of $f$ and $\lambda_c$ is defined by
\bqs
\lambda_c:=\frac{\eta(\mathscr{R}_0-1)}{k+1-2\sqrt{k}}>0.
\eqs
\end{thm}

We remark that there is a threshold on the parameter $\lambda$, which depends on the degree $k$ of the homogeneous tree and the basic reproduction number $\mathscr{R}_0$, below which the stationary solution asymptotically converges to $\mathcal{I}_*>0$ the unique positive zero of $f$ and above which it asymptotically converges to zero. Let us heuristically explain how this critical value emerges. For $n\geq2$, stationary solutions satisfy
\bqs
0=f( \mathcal{I}_n^\infty)+\lambda\left(\mathcal{I}_{n-1}^\infty-(k+1)\mathcal{I}_n^\infty+k\mathcal{I}_{n+1}^\infty\right),
\eqs
and one can look for exponential supersolutions of the form $\overline{\mathcal{I}}_n^\infty=Ce^{-\gamma n}$ for some well chosen constant $C>0$ and $\gamma>0$. Using the concavity of $f$, we get that
\bqs
f( \overline{\mathcal{I}}_n^\infty)+\lambda\left(\overline{\mathcal{I}}_{n-1}^\infty-(k+1)\overline{\mathcal{I}}_n^\infty+k\overline{\mathcal{I}}_{n+1}^\infty\right)< \left( \eta(\mathscr{R}_0-1)+\lambda\left( e^{\gamma}-(k+1)+ke^{-\gamma}\right) \right) Ce^{-\gamma n},
\eqs
and we denote
\bqs
\mathcal{D}(\gamma):=\eta(\mathscr{R}_0-1)+\lambda\left( e^{\gamma}-(k+1)+ke^{-\gamma}\right).
\eqs
We readily remark that $\mathcal{D}$ is convex with $\mathcal{D}(0)=\eta(\mathscr{R}_0-1)>0$ and $\mathcal{D}(+\infty)=+\infty$. Furthermore, $\mathcal{D}'(\gamma)=0$ if and only if $\gamma=\ln \sqrt{k}$, and we note that
\bqs
\mathcal{D}(\ln\sqrt{k})= \eta(\mathscr{R}_0-1)-\lambda(k+1-2\sqrt{k}) = (k+1-2\sqrt{k}) (\lambda_c-\lambda).
\eqs
As a consequence, when $\lambda>\lambda_c$ one can find $\gamma>0$ such that $\mathcal{D}(\gamma)<0$. This implies that
\bqs
f( \overline{\mathcal{I}}_n^\infty)+\lambda\left(\overline{\mathcal{I}}_{n-1}^\infty-(k+1)\overline{\mathcal{I}}_n^\infty+k\overline{\mathcal{I}}_{n+1}^\infty\right)< 0,
\eqs
and $\overline{\mathcal{I}}_n^\infty=Ce^{-\gamma n}$ is a supersolution for $n\geq 2$. For values of $\lambda$ below the critical value $\lambda_c$ the same ingredient as in the proof of Theorem~\ref{thm1} applies and it relies on the so called {\em hair-trigger effect} which holds true for the Fisher-KPP equation in this setting \cite{HH19}.

The dichotomy presented in Theorem~\ref{thm6} greatly differs from the case $k=1$ on the lattice, and is somehow counterintuitive. Indeed, although we are in the case $\mathscr{R}_0>1$, if the intensity of exchanges $\lambda$ is too large, then the epidemic is no longer able to spread into the graph. Actually, as $\lambda_c \sim \frac{\eta (\mathscr{R}_0-1)}{k}\rightarrow 0$ as $k\rightarrow+\infty$, we get that the higher the degree of the tree is the less likely is an epidemic to spread into the network. This paradoxical behavior can be once again explained by noticing that exchange terms in the tree are given by the linear superposition of a discrete diffusive part $\mathcal{I}_{n-1}(t)-2\mathcal{I}_n(t)+\mathcal{I}_{n+1}(t)$ and an advection term $(k-1)(\mathcal{I}_{n+1}(t)-\mathcal{I}_n(t))$ which transports individuals up the root of the tree. For $k$ large, it should then be expected that this advection term dominates the reaction-diffusion part and prevents an epidemic to spread.

We also get that the unique stationary solution is a global attractor for the dynamics \eqref{KPPtree}.

\begin{thm}\label{thm7}
Let $\left(\mathcal{I}_n(t)\right)_{n\geq 1}$ be the solution of \eqref{KPPlike} starting from some nonnegative bounded compactly supported initial condition.
Then  $\left(\mathcal{I}_n(t)\right)_{n\geq1}$ converges as $t\rightarrow+\infty$ locally uniformly to $\left(\mathcal{I}_n^\infty\right)_{n\geq1}$.
\end{thm}

\begin{figure}
\centering
\begin{tabular}{ccc}
\setlength{\tabcolsep}{0pt}
\includegraphics[width=0.32\textwidth]{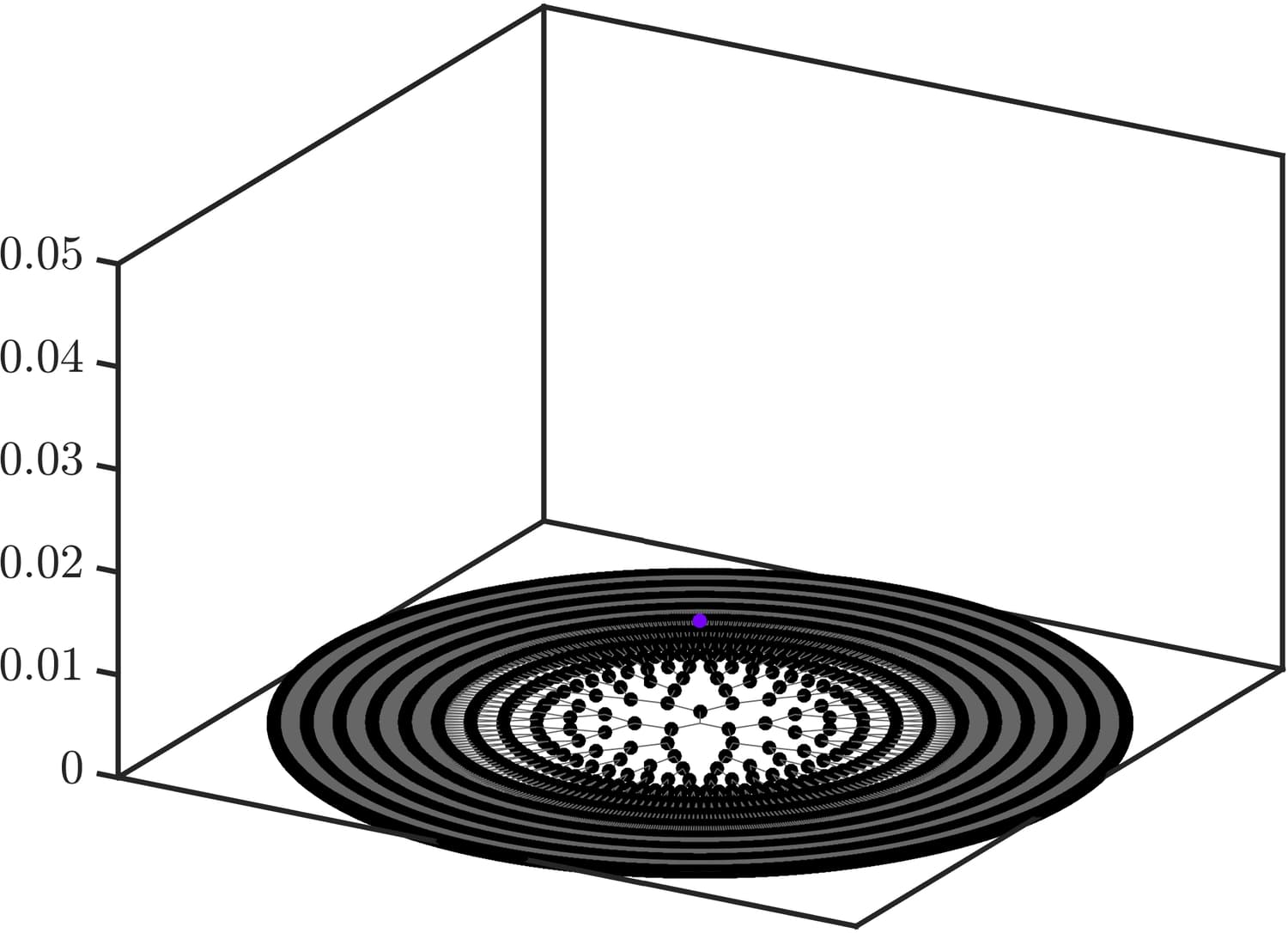}&
\includegraphics[width=0.32\textwidth]{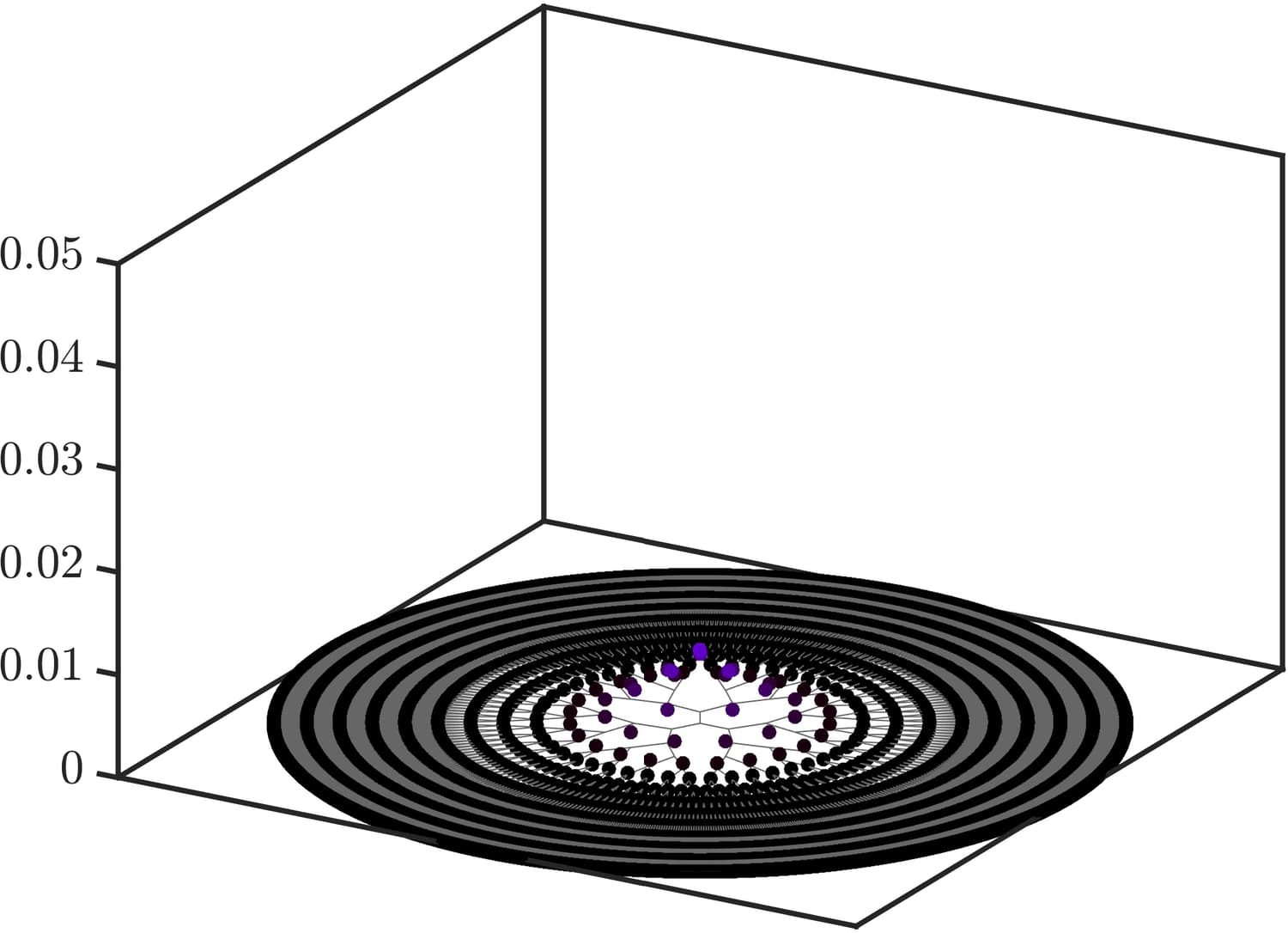}&
\includegraphics[width=0.32\textwidth]{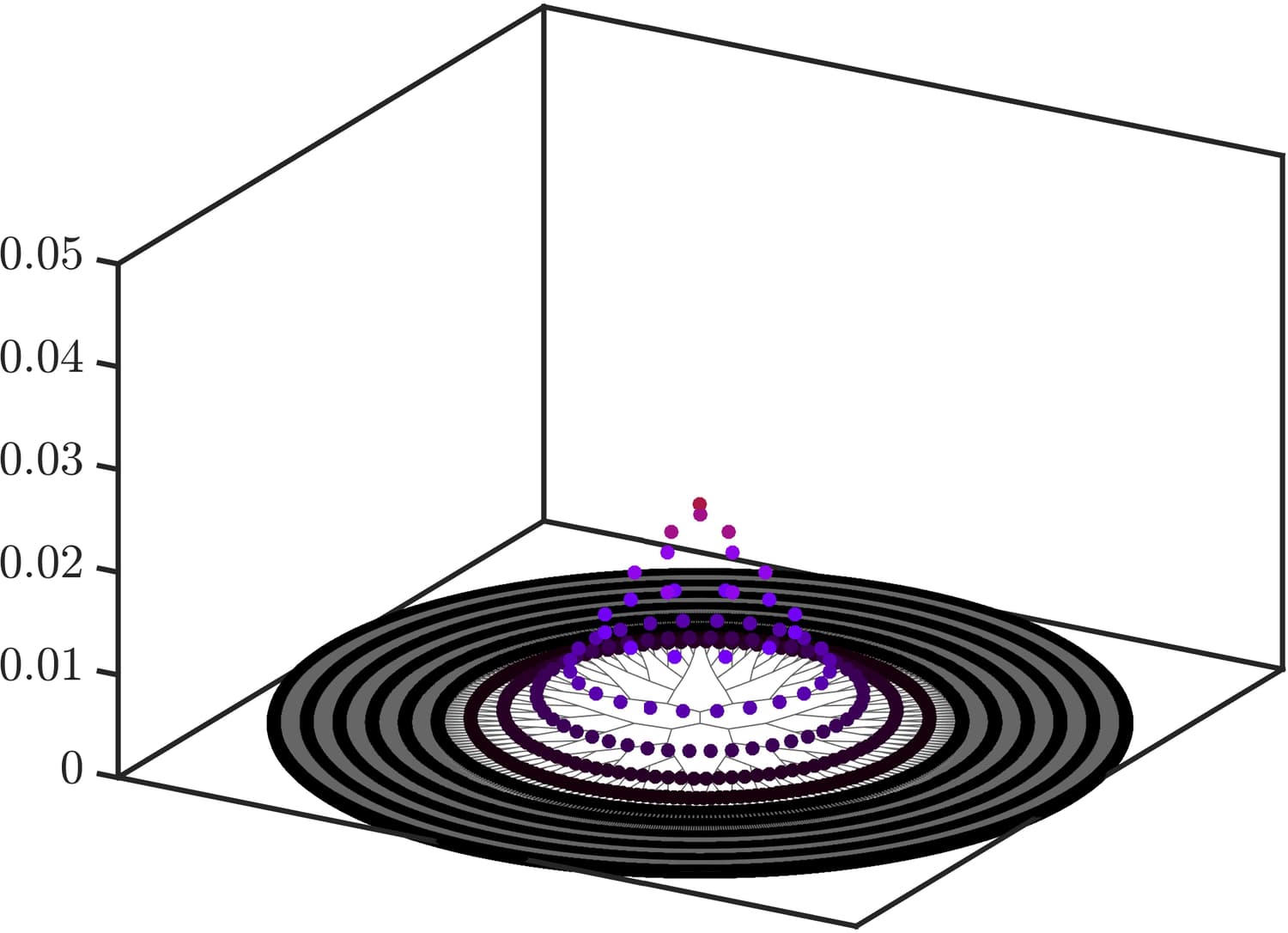}\\
$t=0$ & $t=10$ & $t=20$ \\
\includegraphics[width=0.32\textwidth]{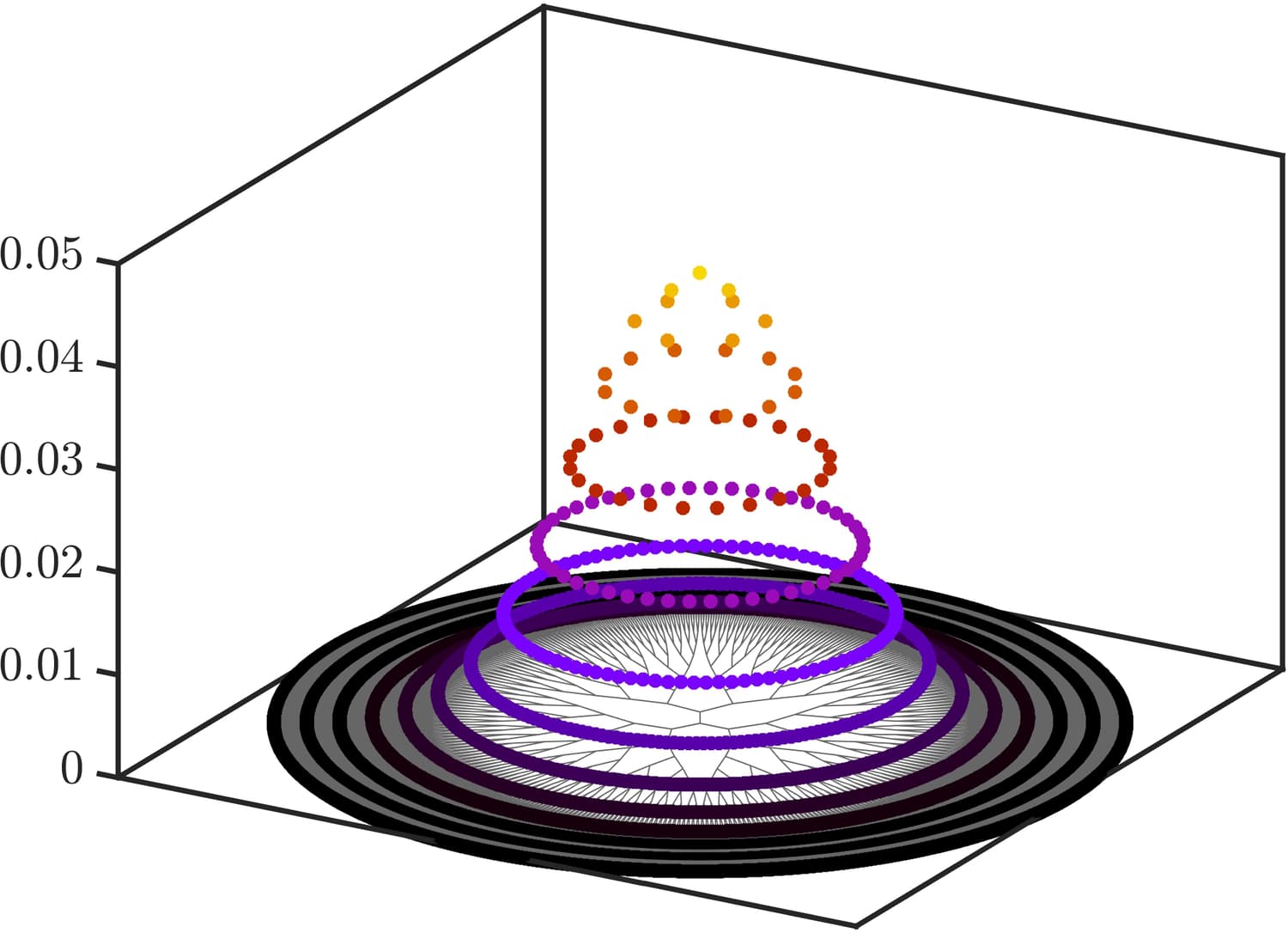}&
\includegraphics[width=0.32\textwidth]{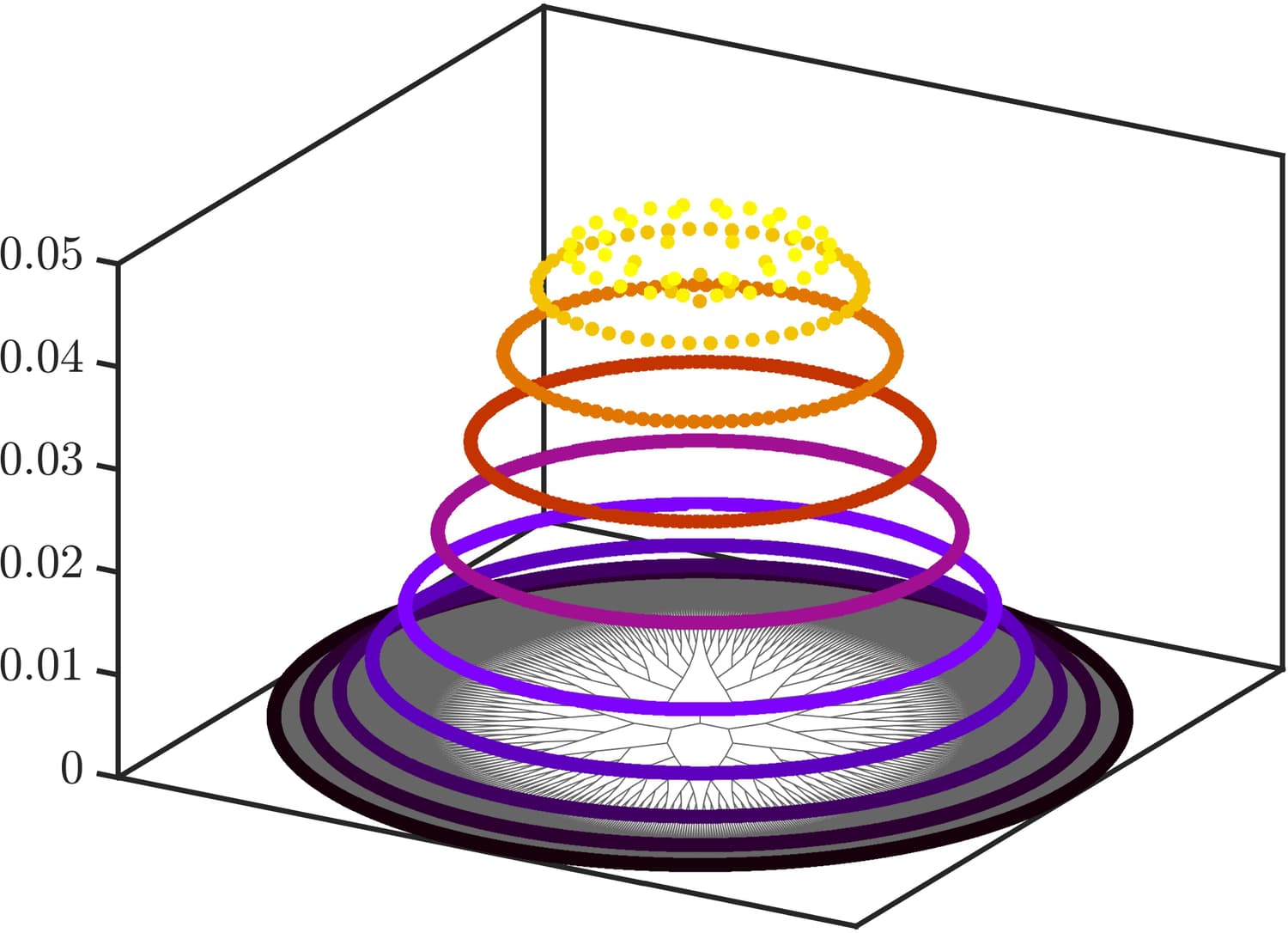}&
\includegraphics[width=0.32\textwidth]{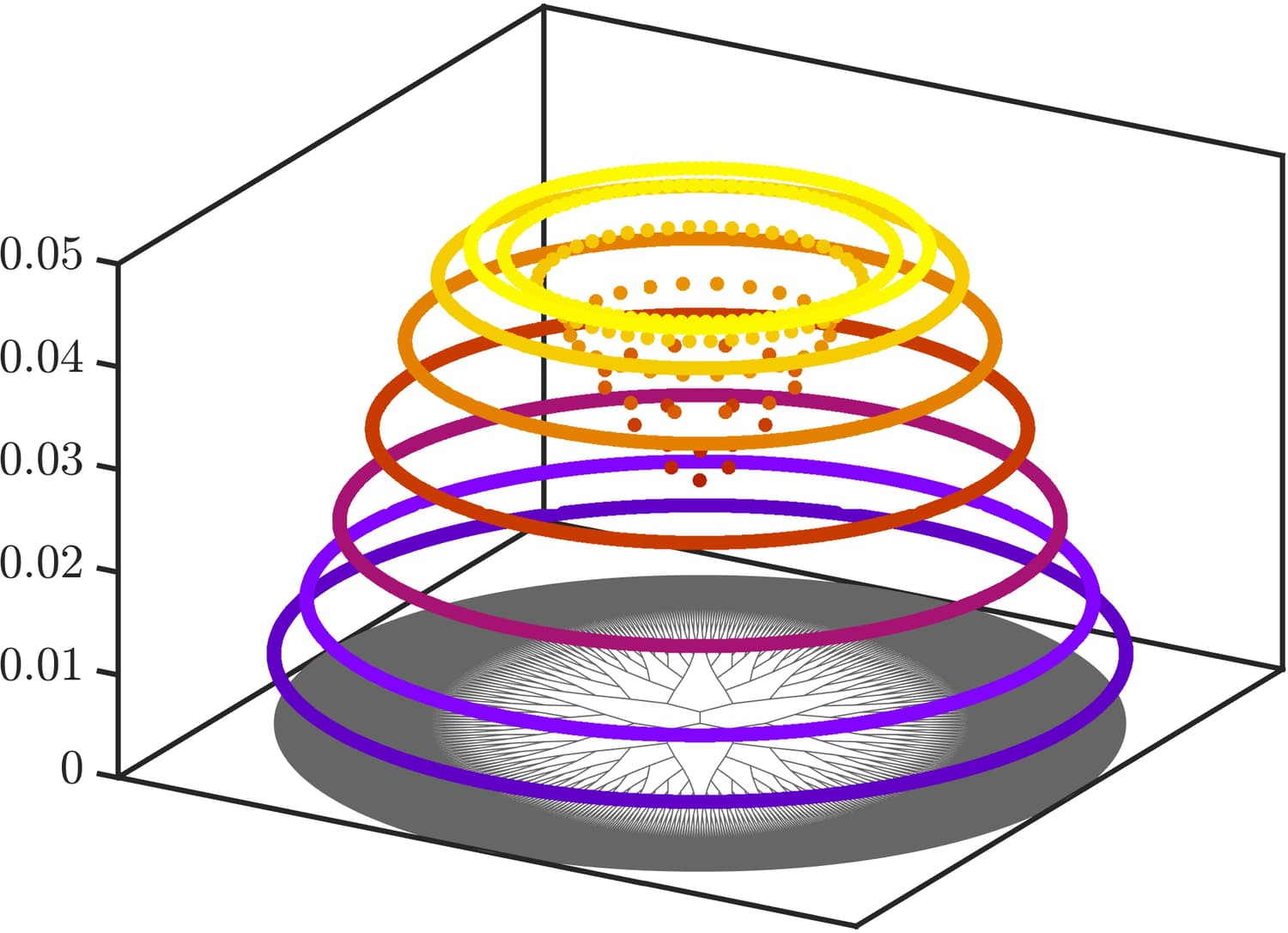}\\
$t=30$ & $t=40$ & $t=50$ \\
\includegraphics[width=0.32\textwidth]{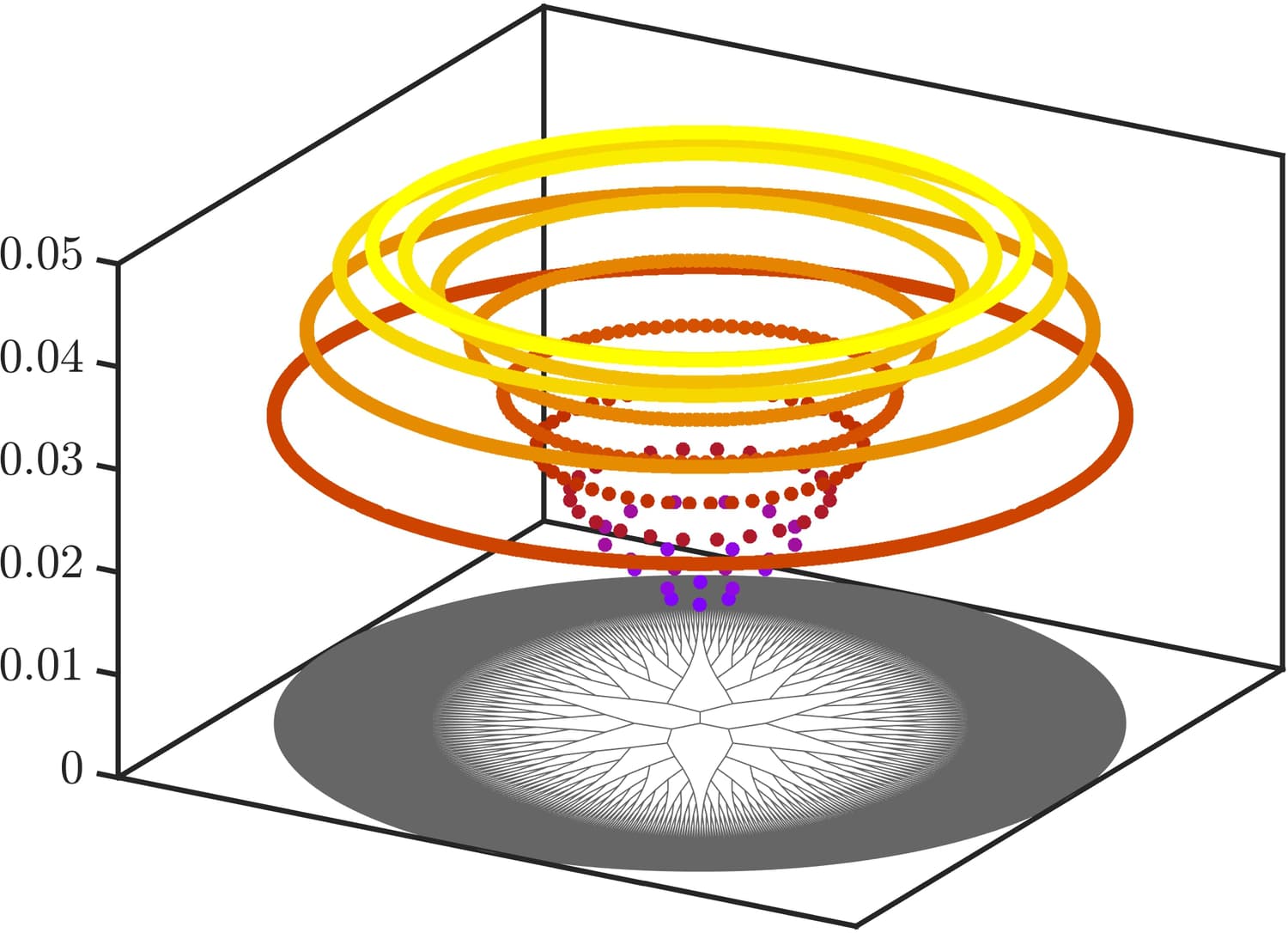}&
\includegraphics[width=0.32\textwidth]{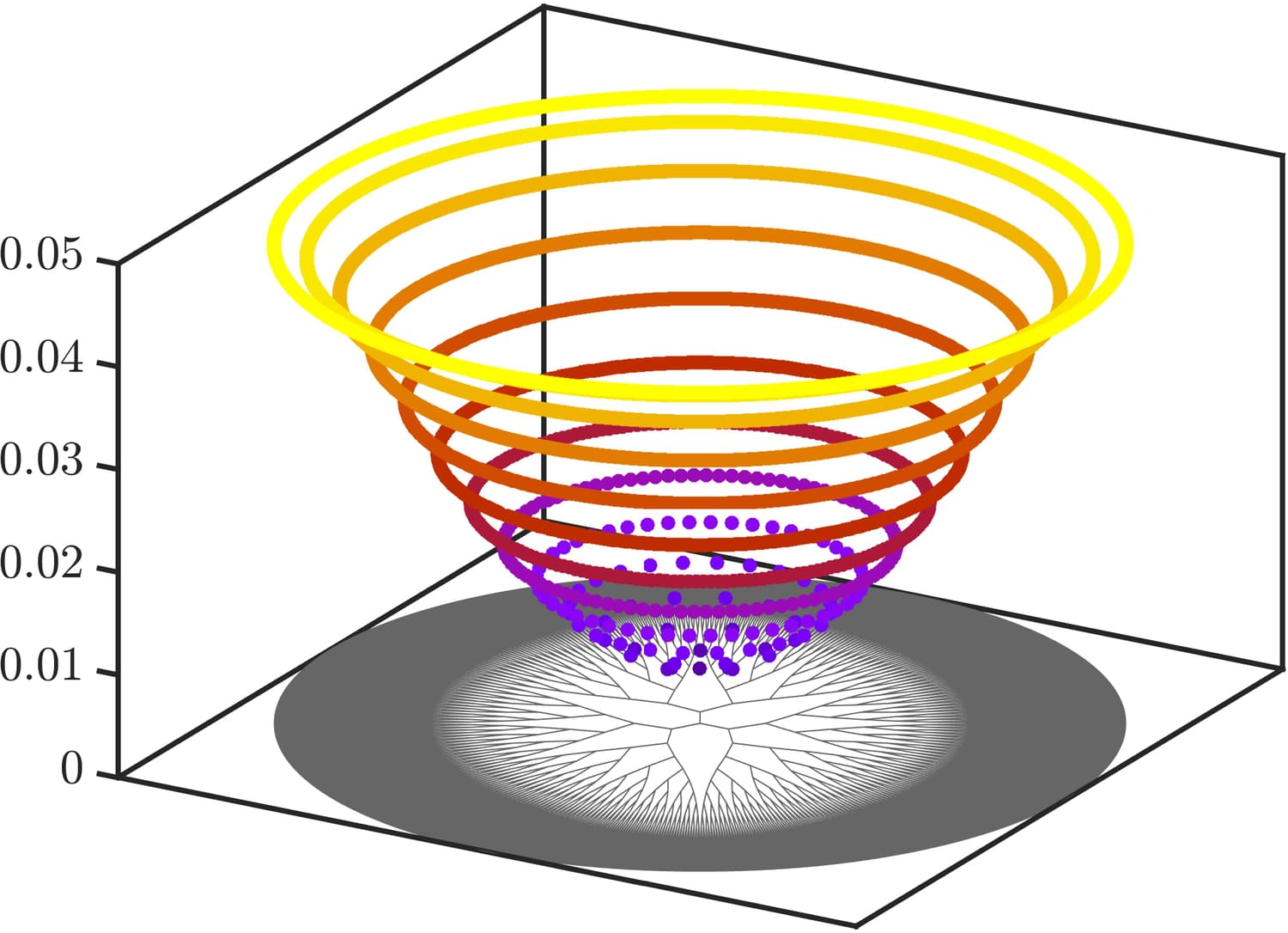}&
\includegraphics[width=0.32\textwidth]{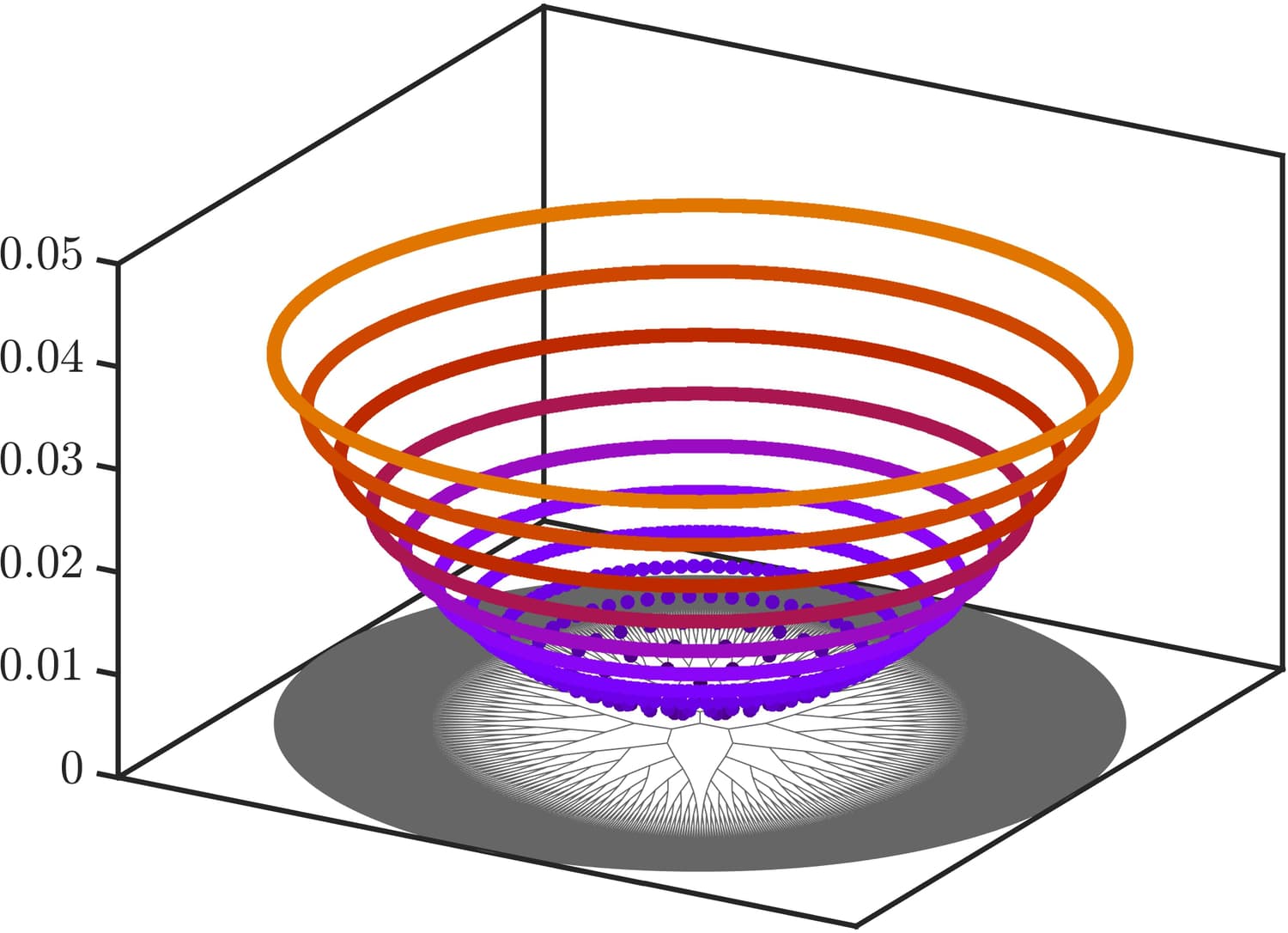}\\
$t=60$ & $t=70$ & $t=80$ \\
\includegraphics[width=0.32\textwidth]{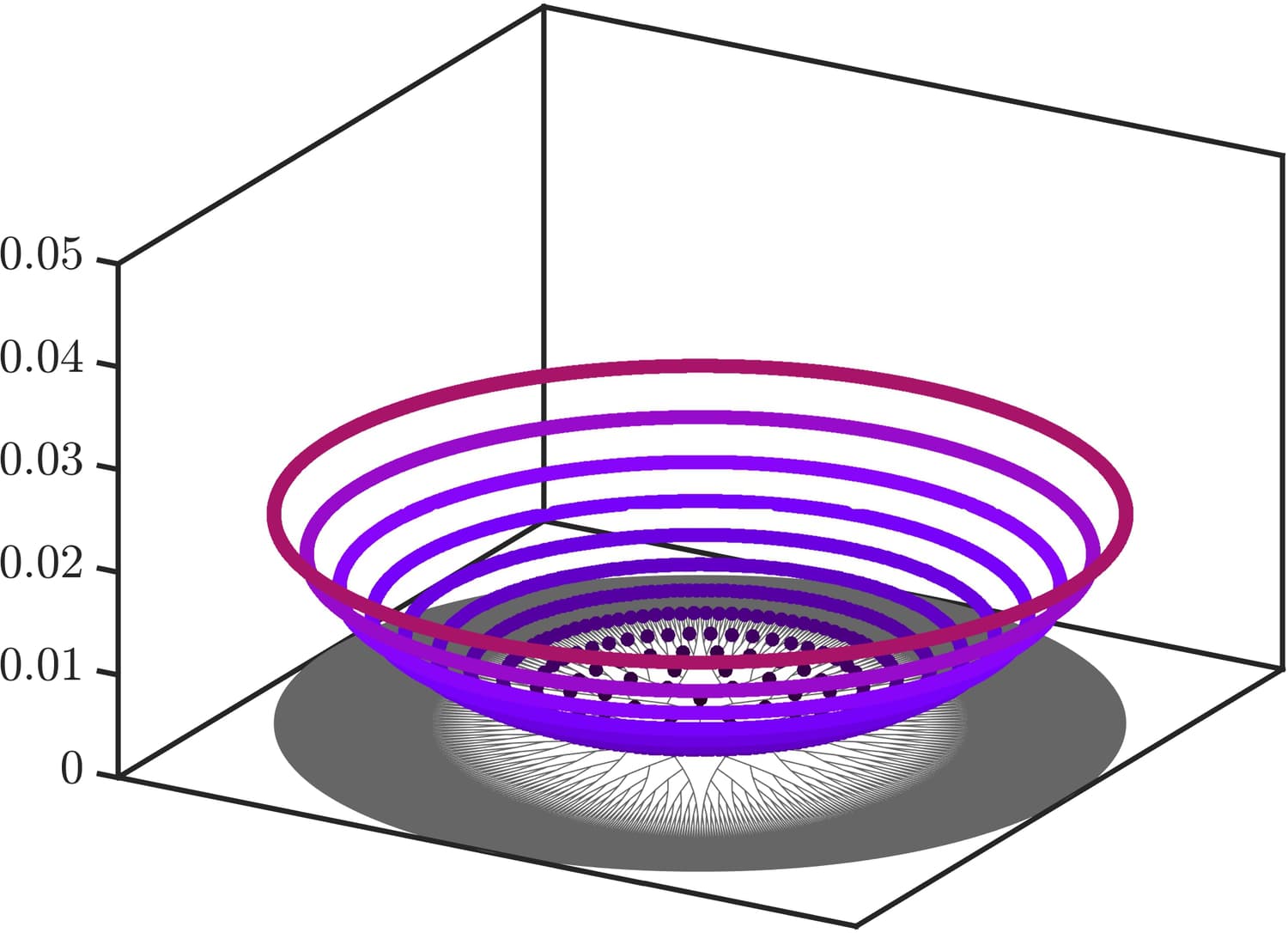}&
\includegraphics[width=0.32\textwidth]{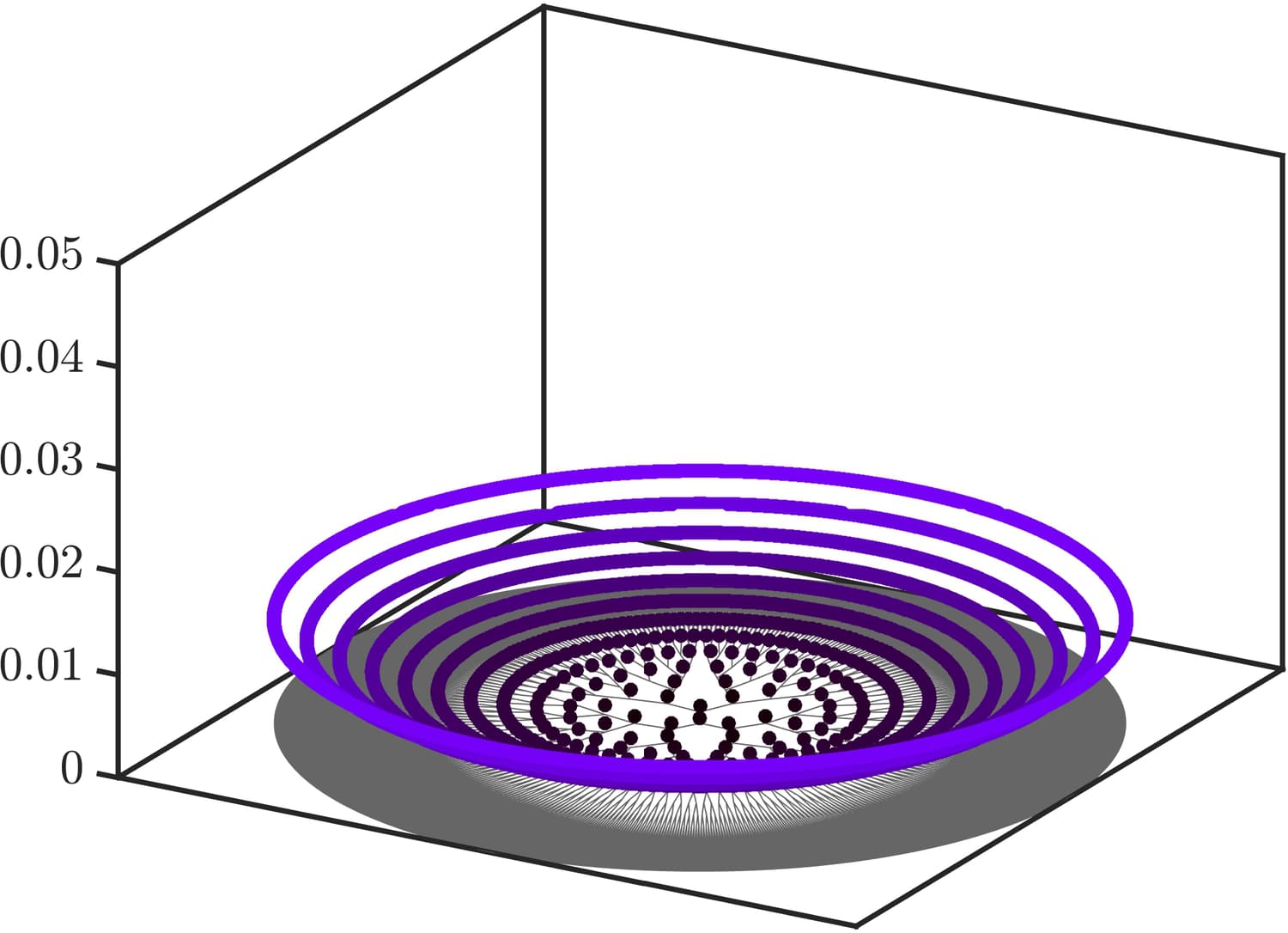}&
\includegraphics[width=0.32\textwidth]{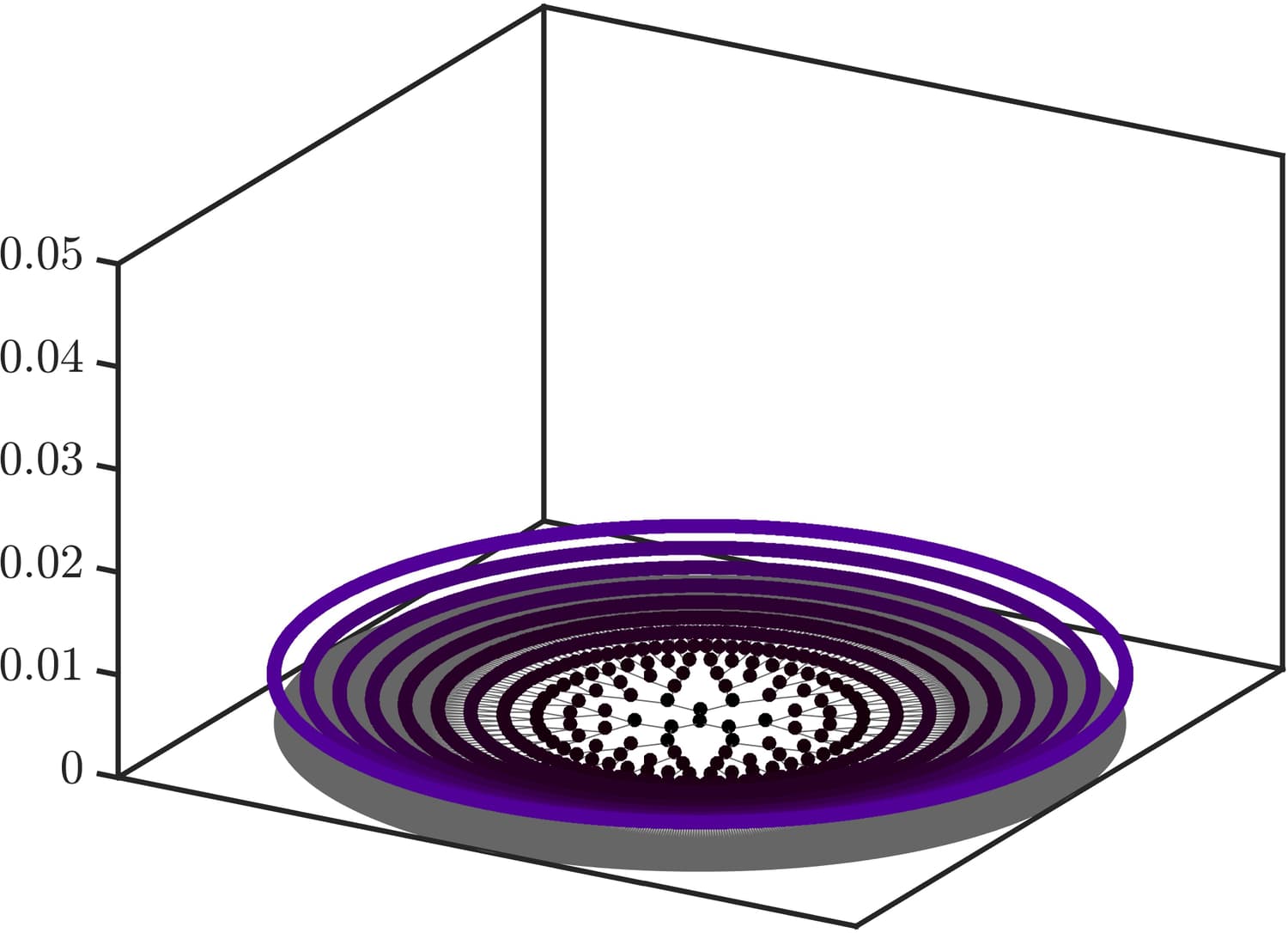}\\
$t=90$ & $t=100$ & $t=110$ 
\end{tabular}
\caption{Time evolution of the density of infected individuals $(I_n(t))_{n\geq1}$ solution of system \eqref{KPPtree} from the initial condition depicted in Figure~\ref{fig:InitialCond} from time $t=0$ to $t=110$ in the case where $\mathscr{R}_0>1$ and for $\lambda\in(0,\lambda_c)$ with $k=2$. We observe the propagation of the infected individuals across the tree. }
\label{fig:Invasion}
\end{figure}

As for the lattice, a consequence of the above Theorem~\ref{thm7} is that the convergence also holds true for the time derivative of $\left(\mathcal{I}_n(t)\right)_{n\geq1}$ such that we can deduce that locally uniformly in $n$ we have
\bqs
I_n(t)=\mathcal{I}_n'(t) \rightarrow 0 \text{ as } t\rightarrow +\infty,
\eqs
which means that the density of infected individuals asymptotically vanishes in time at each fixed lattice site. We can also quantify the density of individuals that will be infected during the course of the epidemic at a given node $n$, that we denote $\mathcal{I}_n^{tot}$. It is given by
\bqs
\mathcal{I}_n^{tot}=s_0\left(1-e^{-\tau \mathcal{I}_n^\infty}\right), \quad n\geq1.
\eqs
Now, using the results of Theorem~\ref{thm5} and Theorem~\ref{thm6}, we obtain the following characterization
\bqs
\underset{|j|\rightarrow+\infty}{\lim}\mathcal{I}_j^{tot}=\left\{
\begin{array}{lcl}
0, & \text{if}  & \mathscr{R}_0 \leq 1, \\ 
s_0\left(1-e^{-\tau \mathcal{I}_*}\right), & \text{if} & \mathscr{R}_0>1 \text{ and } 0<\lambda <\lambda_c, \\
0,& \text{if} & \mathscr{R}_0>1 \text{ and } \lambda >\lambda_c.
\end{array}
\right.
\eqs
As a consequence, an epidemic can only spread in homogeneous trees of degree $k\geq2$ when $\mathscr{R}_0>1$ and $0<\lambda <\lambda_c$. We present in Figure~\ref{fig:Invasion} the time evolution of the density of infected individuals $(I_n(t))_{n\geq1}$ solution of system \eqref{KPPtree} from the initial condition depicted in Figure~\ref{fig:InitialCond} from time $t=0$ to $t=110$ in the case where $\mathscr{R}_0>1$ and for $\lambda\in(0,\lambda_c)$ with $k=2$. We observe the propagation of the infected individuals across the tree.

In that case, we can also characterize at which speed the epidemic spreads into the tree. 

\begin{thm}\label{thm8}
Assume that $\mathscr{R}_0>1$ and $0<\lambda<\lambda_c$. We define $c_*^k>0$ as
\bqs
c_*^k:=\underset{\gamma>0}{\min}~  \frac{\eta \left( \mathscr{R}_0-1\right)+\lambda \left( e^{\gamma}-(k+1)+ke^{-\gamma} \right)}{\gamma}.
\eqs
Then, the solution $\left(\mathcal{I}_n(t)\right)_{n\geq1}$ of \eqref{KPPtree}-\eqref{KPPtreeIC} satisfies:
\begin{itemize}
\item[(i)] $\forall c\in(0,c_*^k)$,
\bqs
\underset{t\rightarrow+\infty}{\lim} \left( \underset{1 \leq n \leq ct}{\sup}\left| \mathcal{I}_n(t)-\mathcal{I}_n^\infty \right| \right) =0;
\eqs
\item[(ii)] $\forall c>c_*^k$,
\bqs
\underset{t\rightarrow+\infty}{\lim} \left( \underset{n \geq ct}{\sup}\left| \mathcal{I}_n(t) \right| \right) =0.
\eqs
\end{itemize}
\end{thm}

\begin{figure}[t!]
  \centering
  \includegraphics[width=.5\textwidth]{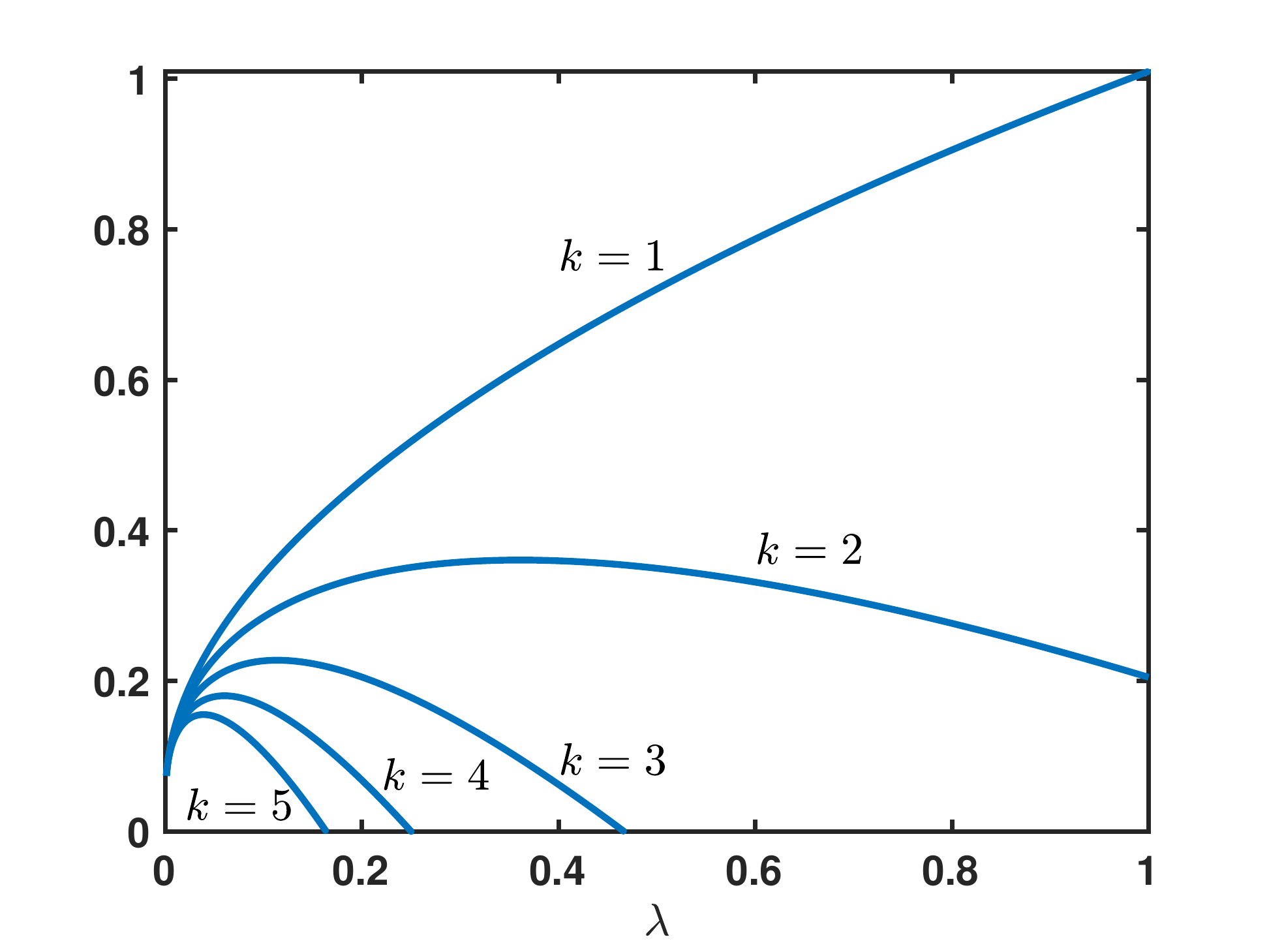} 
  \caption{Linear spreading speed $c_*^k$ of a homogeneous tree $\mathbb{T}_k$ of degree $k$ for $k\in\left\{1,\cdots,5\right\}$ as a function of $\lambda$ when all other parameters are fixed such that $\mathscr{R}_0=\frac{\tau S_0}{\eta}>1$. When $k=1$, that is on the lattice, the spreading speed is a monotone function of $\lambda$. On the other, for $k\geq 2$ the spreading speed is not monotone and there exists some critical value $\lambda_c$ for which the spreading speed vanishes.}
  \label{fig:clintree}
\end{figure}

We have presented in Figure~\ref{fig:clintree} the spreading speed $c_*^k$ as a function of the parameter $\lambda$ for several values of $k$, with the convention that when $k=1$, we have $c_*^1=c_*$ where $c_*$ is given by Theorem~\ref{thm3}. We recover that when $k=1$, that is on the lattice, the spreading speed is a monotone function of $\lambda$. On the other hand, when $k\geq2$, we observe two key features. First, the spreading speed is no longer monotone. It is strictly increasing up to some critical value $\lambda_0>0$ and then strictly decreasing. Second, the spreading speed $c_*^k$ vanishes at the some critical value of the parameter $\lambda$, and it turns out that this critical value is precisely $\lambda_c$ from Theorem~\ref{thm6}. We summarize these properties in the following proposition whose proof can be found in \cite{HH19}. 

\begin{prop}
Assume that $\mathscr{R}_0>1$. Let $c_*^k$ be the spreading speed defined in Theorem~\ref{thm8}. Then we have:
\begin{itemize}
\item $c_*^k>0$ for all $\lambda \in (0,\lambda_c)$ with $c_*^k=0$ when $\lambda=\lambda_c$;
\item $\lambda\mapsto c_*^k$ is monotone increasing on $(0,\lambda_0)$ and monotone decreasing on $(\lambda_0,\lambda_c)$ where 
\bqs
\lambda_0:=\frac{\eta(\mathscr{R}_0-1)}{(k-1)\ln k}.
\eqs
\end{itemize}
\end{prop}

The above result shows that $\lambda_0$ maximizes the speed at which the epidemic spreads into the tree. At a given fixed basic reproduction number $\mathscr{R}_0>1$, this critical value decreases as $k$ increases, and the associated spreading speed is
\bqs
c_*^k(\lambda=\lambda_0) = \frac{\eta(\mathscr{R}_0-1)}{\ln k}.
\eqs 
This is somehow counterintuitive as one would have expected that the spreading speed would increase with the degree $k$ of the tree.

\subsection{Outline of the paper}

The remainder of the paper is organized as follows. In Section~\ref{secZ}, we present the proofs of our main results regarding the case $\G=\Z$. In the following Section~\ref{secT}, we explain how these proofs transpose to the case of a homogeneous tree $\G=\mathbb{T}_k$ for $k\geq2$. We conclude with a discussion in Section~\ref{disc}.

\section{Spreading properties on the lattice $\Z$}\label{secZ}

In this section, we study the spreading properties of \eqref{KPPlike} set on the lattice $\Z$. Our approach is similar to the continuous case and we refer for example to \cite{DG14,BRR21}. It relies on comparison principles techniques and some known results for the Fisher-KPP equation on the lattice \cite{HH19}.

\subsection{Proof of Theorem~\ref{thm1}}

In section, we study the existence of stationary solutions to \eqref{KPPlike}, that is we look for bounded sequences $\left(\mathcal{I}_j\right)_{j\in\Z} \in \ell_\infty(\Z)$ that satisfy
\bqq
0=f(\mathcal{I}_j)+I_j^0+\lambda\left( \mathcal{I}_{j-1}-2\mathcal{I}_j+\mathcal{I}_{j+1}\right),
\label{stat}
\eqq
for each $j\in\Z$. We first remark that the $0$ sequence is a subsolution to the above equation since we assume that $I_j^0\geq0$. On the other hand has $f(+\infty)=-\infty$ and each $I_j^0$ is bounded, there exists some positive real $\rho>0$ such that the constant sequence with $\mathcal{I}_j=\rho$ for each $j\in\Z$ is a supersolution of \eqref{stat}. We denote by $\left(\overline{\mathcal{I}_j}(t)\right)_{j\in\Z}$ the time dependent solution of the Cauchy problem
\bqq
\left\{
\begin{split}
\mathcal{I}_j'(t)&=f(\mathcal{I}_j(t))+I_j^0+\lambda\left( \mathcal{I}_{j-1}(t)-2\mathcal{I}_j(t)+\mathcal{I}_{j+1}(t)\right),\\
\mathcal{I}_j(0)&= \rho,\\
\end{split}
\right.
\label{Cauchy}
\eqq
for each $j\in\Z$ and $t>0$. Since the constant sequence with $\mathcal{I}_j=\rho$ for each $j\in\Z$ is a supersolution of \eqref{stat} it is also a supersolution of \eqref{Cauchy}, and it follows from the comparison principle that $\left(\overline{\mathcal{I}_j}(t)\right)_{j\in\Z}$ is nonincreasing in the time variable. Furthermore, it satisfies $0\leq \overline{\mathcal{I}_j}(t) \leq \rho$ for each $j\in\Z$ and $t>0$. Therefore $\left(\overline{\mathcal{I}_j}(t)\right)_{j\in\Z}$ converges as $t\rightarrow+\infty$ to some sequence $\left(\mathcal{I}_j^\infty\right)_{j\in\Z}$. Thanks to the regularity in time of the solution of the Cauchy problem \eqref{Cauchy}, this stationary sequence is solution of \eqref{stat} and by construction $0\leq \mathcal{I}_j^\infty \leq \rho$ for each $j\in\Z$. As $I_j^0>0$ for some $j$, we have that the $0$ sequence is not a solution of \eqref{stat}. We claim that this implies that $0< \mathcal{I}_j^\infty$ for each $j\in\Z$. Indeed, assume by contradiction that there exists some $j_0\in\Z$ for which $\mathcal{I}_{j_0}^\infty=0$, the equation \eqref{stat} gives
\bqs
0 = I_{j_0}^0+\lambda\left(\mathcal{I}_{j_0-1}^\infty+\mathcal{I}_{j_0+1}^\infty \right).
\eqs
If $j_0$ belongs to the support of $\left( I^0_j \right)_{j\in\Z}$ then we have reached a contradiction. Otherwise, we deduce that $\mathcal{I}_{j_0-1}^\infty=\mathcal{I}_{j_0+1}^\infty=0$. We can then repeat the argument such that $j_0+k$ belongs to the support of $\left( I^0_j \right)_{j\in\Z}$ which then leads to a contradiction and proves the claim.

We now derive the limit as $|j|\rightarrow+\infty$. We introduce the sequence of shifts $\mathcal{I}_{j+\tau_n}^\infty$ with $\tau_n \in \Z$ for $n\in\N$ and $|\tau_n|\rightarrow+\infty$ as $n\rightarrow+\infty$. Up to subsequences, it converges towards a nonnegative bounded sequence $\left(\widetilde{\mathcal{I}}_j\right)_{j\in\Z}$ which satisfies 
\bqs
0=f(\widetilde{\mathcal{I}}_j)+\lambda\left( \widetilde{\mathcal{I}}_{j-1}-2\widetilde{\mathcal{I}}_j+\widetilde{\mathcal{I}}_{j+1}\right), \quad j\in\Z.
\eqs
Next, we remark that if $\mathscr{R}_0\leq1$, then $f'(0)\leq 0$ and $f<0$ on $(0,+\infty)$ from which we readily deduce that $\widetilde{\mathcal{I}}_j=0$ for all $j\in\Z$. This proves that $\underset{|j|\rightarrow+\infty}{\lim}\mathcal{I}_j^\infty=0$ when $\mathscr{R}_0\leq1$. Now, if $\mathscr{R}_0>1$ we note that the sequence $\left(\mathcal{I}_j^\infty\right)_{j\in\Z}$ is a supersolution to the Fisher-KPP equation set on the lattice $\Z$
\bqs
\mathcal{I}_j'(t)=f(\mathcal{I}_j(t))+\lambda\left( \mathcal{I}_{j-1}(t)-2\mathcal{I}_j(t)+\mathcal{I}_{j+1}(t)\right), \quad j\in\Z, \quad t>0.
\eqs
For this equation, we know that any solution with a positive, bounded initial condition converges locally uniformly on $\Z$ as $t\rightarrow+\infty$ to the positive zero of $f$ which is denoted $\mathcal{I}_*>0$ (see \cite{W82,HH19}). By comparison, we must have $\mathcal{I}_j^\infty \geq \mathcal{I}_*$ for each $j\in\Z$. This also shows that necessarily $\widetilde{\mathcal{I}}_j \geq \mathcal{I}_*$ and we infer that $\widetilde{\mathcal{I}}_j = \mathcal{I}_*$ since $f<0$ on $(\mathcal{I}_*,+\infty)$. As a consequence, we have proved that $\underset{|j|\rightarrow+\infty}{\lim}\mathcal{I}_j^\infty=\mathcal{I}_*$ when $\mathscr{R}_0>1$.

To conclude the proof of Theorem~\ref{thm1}, it remains to show that \eqref{stat} has a unique stationary solution. We distinguish between the cases $\mathscr{R}_0>1$ and $\mathscr{R}_0\leq1$.

\begin{itemize}
\item Case $\mathscr{R}_0>1$. Let $\left(\mathcal{I}_j^\infty\right)_{j\in\Z}$ and $\left(\mathcal{J}_j^\infty\right)_{j\in\Z}$ be two positive, bounded solutions to \eqref{stat}. Since $\mathcal{I}_j^\infty,\mathcal{J}_j^\infty \geq \mathcal{I}_*$ and both $\mathcal{I}_j^\infty$ and $\mathcal{J}_j^\infty$ are bounded for each $j$, we have that
\bqs
\theta:=\underset{j\in\Z}{\sup} ~ \frac{\mathcal{I}_j^\infty}{\mathcal{J}_j^\infty}>0,
\eqs
is a well-defined quantity. Assume by contradiction that $\theta>1$. Since $\underset{|j|\rightarrow+\infty}{\lim}\mathcal{I}_j^\infty=\underset{|j|\rightarrow+\infty}{\lim}\mathcal{J}_j^\infty=\mathcal{I}_*$, we have that the above supremum is a maximum attained at some $j_0\in\Z$. We obtain by subtracting the equations evaluated at $j=j_0$ that  
\bqs
0=f(\mathcal{I}_{j_0}^\infty)-\theta f(\mathcal{J}_{j_0}^\infty) +\underbrace{I_{j_0}^0(1-\theta)}_{<0}+\underbrace{\lambda\left( \mathcal{I}_{{j_0}-1}^\infty-2\mathcal{I}_{j_0}^\infty+\mathcal{I}_{{j_0}+1}^\infty\right)-\theta \lambda\left( \mathcal{J}_{{j_0}-1}^\infty-2\mathcal{J}_{j_0}^\infty+\mathcal{J}_{{j_0}+1}^\infty\right)}_{\leq 0}.
\eqs
This yields that
\bqs
\theta f(\mathcal{J}_{j_0}^\infty) \leq f(\mathcal{I}_{j_0}^\infty) = f(\theta \mathcal{J}_{j_0}^\infty),
\eqs
which is impossible by concavity of $f$. As a consequence, we have proved that $\mathcal{I}_j^\infty\leq \mathcal{J}_j^\infty$ for each $j\in\Z$. By reversing the role of the two solutions, we reach the conclusion.
\item Case $\mathscr{R}_0\leq1$. Let $\left(\mathcal{I}_j^\infty\right)_{j\in\Z}$ and $\left(\mathcal{J}_j^\infty\right)_{j\in\Z}$ be two positive, bounded solutions to \eqref{stat}. Note that they both tend to $0$ as $|j|\rightarrow +\infty$. Take $\epsilon >0$ and denote $\left(\mathcal{J}_j^{\infty,\epsilon}\right)_{j\in\Z}$ the sequence with $\mathcal{J}_j^{\infty,\epsilon}=\mathcal{J}_j^{\infty}+\epsilon$ which satisfies
\bqs
 f(\mathcal{J}_j^{\infty,\epsilon})+I^0_j+\lambda\left(\mathcal{J}_{j-1}^{\infty,\epsilon}-2\mathcal{J}_j^{\infty,\epsilon}+\mathcal{J}_{j+1}^{\infty,\epsilon}\right) <0, \quad j\in\Z,
\eqs
since $f$ is decreasing on $\R_+$. Assuming by contradiction that $\mathcal{I}_j^{\infty}>\mathcal{J}_j^{\infty,\epsilon}$ for some $j$, and subtracting the  equation for $\mathcal{I}_j^{\infty}$ and the strict inequality for $\mathcal{J}_j^{\infty,\epsilon}$, we end up with the inequality
\bqs
\underbrace{f(\mathcal{I}_{j_0}^\infty)- f(\mathcal{J}_{j_0}^{\infty,\epsilon})}_{<0} +\underbrace{\lambda\left( \mathcal{I}_{{j_0}-1}^\infty-2\mathcal{I}_{j_0}^\infty+\mathcal{I}_{{j_0}+1}^\infty\right)- \lambda\left( \mathcal{J}_{{j_0}-1}^{\infty,\epsilon}-2\mathcal{J}_{j_0}^{\infty,\epsilon}+\mathcal{J}_{{j_0}+1}^{\infty,\epsilon}\right)}_{\leq 0}>0
\eqs
for some $j_0\in\Z$ where the maximum is attained. And once again, we have reached a contradiction, thus $\mathcal{I}_j^{\infty} \leq \mathcal{J}_j^{\infty,\epsilon}$ for all $j\in\Z$ and arbitrary $\epsilon>0$ which gives $\mathcal{I}_j^{\infty} \leq \mathcal{J}_j^{\infty}$ for all $j\in\Z$ and concludes the proof in that case.
\end{itemize}

\subsection{Proof of Theorem~\ref{thm2}}

Let $\left(\mathcal{I}_j(t)\right)_{j\in\Z}$ be the solution of \eqref{KPPlike} starting respectively from some nonnegative bounded compactly supported initial condition $\left(\mathcal{I}_j^0\right)_{j\in\Z}$. Then, we denote by $\left(\underline{\mathcal{I}}_j(t)\right)_{j\in\Z}$ and $\left(\overline{\mathcal{I}_j}(t)\right)_{j\in\Z}$ the time dependent solutions of \eqref{KPPlike} starting with initial condition being the $0$ constant sequence and the constant sequence
 with $\overline{\mathcal{I}_j}(0)=\rho$ where $\rho>0$ is chosen large enough such that both $f(\rho)+\max I_j^0<0$ and $\rho>\max \mathcal{I}_j^0$. Then we have for all $t>0$
\bqs
0<\underline{\mathcal{I}}_j(t)\leq \mathcal{I}_j(t) \leq \overline{\mathcal{I}_j}(t)<\rho, \quad j\in\Z.
\eqs
By comparison, $\underline{\mathcal{I}}_j(t)$ and $\overline{\mathcal{I}_j}(t)$ are respectively increasing and decreasing in time and converge locally uniformly to two stationary solutions $(\underline{\mathcal{I}}_j^\infty)_{j\in\Z}$ and $(\overline{\mathcal{I}_j}^\infty)_{j\in\Z}$ of \eqref{stat}. By Theorem~\ref{thm1}, we must have $\underline{\mathcal{I}}_j^\infty=\overline{\mathcal{I}}_j^\infty=\mathcal{I}_j^\infty$ where $(\mathcal{I}_j^\infty)_{j\in\Z}$ is the unique stationary solution of \eqref{stat}. The proof is thereby complete.

\subsection{Proof of Theorem~\ref{thm3}}

In this section we assume that $\mathscr{R}_0>1$. And we denote $\left(\mathcal{I}_j(t)\right)_{j\in\Z}$ the solution of \eqref{KPPlike}-\eqref{ICkpp}. We let $c\in(0,c_*)$ and consider a sequence $(t_n)_{n\in\N}$ such that $t_n\rightarrow+\infty$ as $n\rightarrow+\infty$ and a sequence $(j_n)_{n\in\N}$ in $\Z$ such that $|j_n|\leq (c_*-c)t_n$. If $(j_n)_{n\in\N}$ is bounded, we know from Theorem~\ref{thm2} that $\mathcal{I}_{j_n}(t_n)-\mathcal{I}_{j_n}^\infty\rightarrow 0\rightarrow 0$ as $n\rightarrow+\infty$ by local uniform convergence. Suppose that up to subsequences $(j_n)_{n\in\N}$ diverges. We recall that the solution $\left(\mathcal{I}_j(t)\right)_{j\in\Z}$ is a supersolution of the Fisher-KPP equation set on the lattice for which spreading occurs with the asymptotic speed $c_*$. We infer that
\bqs
\underset{n\rightarrow+\infty}{\liminf} \left( \mathcal{I}_{j_n}(t_n)-\mathcal{I}_{j_n}^\infty \right) \geq \mathcal{I}_*-\mathcal{I}_*=0.
\eqs
Next, we let $\gamma_0>0$ be the unique positive zero of  $\varphi(\gamma):=f'(\mathcal{I}_*)+\lambda \left( e^{-\gamma}+e^\gamma-2\right)$. We define the sequence $\mathcal{H}_j:=\mathcal{I}_*+ \beta e^{-\gamma_0 j}$ for $j\in\Z$. We readily remark that outside the support of $(I_j^0)$, we have
\bqs
f(\mathcal{H}_j)+\lambda \left( \mathcal{H}_{j-1}-2\mathcal{H}_j+\mathcal{H}_{j+1}\right)<\left( f'(\mathcal{I}_*)+\lambda \left( e^{-\gamma_0}+e^\gamma_0-2\right) \right)\beta e^{-\gamma_0 j}=0,
\eqs
such that $(\mathcal{H}_j)_{j\in\Z}$ is a supersolution outside the support of $(I_j^0)$ thanks to the concavity of $f$. Then, we select $\beta>0$ large enough, so that $\mathcal{I}_j(t)<\mathcal{H}_j$ for each $j \in \mathrm{supp} (I_j^0)$ and for all $t>0$. Hence, we obtain by comparison that $\mathcal{I}_j(t)<\mathcal{H}_j$ for each $j\in\Z$ and $t>0$. A reflection symmetry implies that, with the same $\beta$, we have $\mathcal{I}_j(t)\leq \mathcal{I}_*+\beta e^{-\gamma_0 |j|}$ for each $j\in\Z$ and $t>0$. It follows that
\bqs
\underset{n\rightarrow+\infty}{\limsup}~ \left( \mathcal{I}_{j_n}(t_n) -\mathcal{I}_{j_n}^\infty\right)\leq \underset{n\rightarrow+\infty}{\limsup}~ \left(\mathcal{I}_*-\mathcal{I}_{j_n}^\infty+\beta e^{-\gamma_0 |j_n|} \right) =0.
\eqs
This concludes the proof of the first statement of the theorem.

We now proceed with the second case. We recall the definition of $c_*>0$ as 
\bqs
c_*:=\underset{\gamma>0}{\min}~  \frac{\eta \left( \mathscr{R}_0-1\right)+\lambda \left( e^{-\gamma}-2+e^{\gamma} \right)}{\gamma},
\eqs
where the minimum is achieved at a unique value $\gamma_*>0$. We introduce the sequence $\mathcal{H}_j(t):=\alpha e^{-\gamma_* (j-c_*t)}$ for $j\in\Z$ for some constant $\alpha>0$ which is fixed large enough such that $\mathcal{I}_j(t)<\mathcal{H}_j(t)$ for each $j \in \mathrm{supp} (I_j^0)$ and for all $t>0$. Outside the support of $(I_j^0)$, we compute
\bqs
\mathcal{H}_j'(t)-f(\mathcal{H}_j(t))-\lambda\left( \mathcal{H}_{j-1}(t)-2\mathcal{H}_j(t)+\mathcal{H}_{j+1}(t)\right)>\left(-c_*\gamma_* -f'(0)-\lambda \left( e^{-\gamma_*}+e^{\gamma_*}-2\right) \right) \mathcal{H}_j(t) =0,
\eqs
thanks to the concavity of $f$. Thus $(\mathcal{H}_j(t))_{j\in\Z}$ is a supersolution outside the support of $(I_j^0)$ and we can deduce that $\mathcal{I}_j(t)<\mathcal{H}_j(t)$ for all $t>0$ and $j\in\Z$. As a consequence, for each $c>c_*$, we have
\bqs
\underset{t\rightarrow+\infty}{\lim} \left( \underset{j \geq ct}{\sup}\left| \mathcal{I}_j(t) \right| \right) =0.
\eqs
A symmetry argument also shows that
\bqs
\underset{t\rightarrow+\infty}{\lim} \left( \underset{j \leq -ct}{\sup}\left| \mathcal{I}_j(t) \right| \right) =0.
\eqs
This concludes the proof of the theorem.

\subsection{Traveling fronts}

We now turn our attention to the existence of traveling wave solutions to system \eqref{TWlattice} such that $0<S(x)<s_0$ and $I(x)>0$ is bounded on $\R$ with asymptotic conditions $S(+\infty)=s_0$ and $I(\pm\infty)=0$. From \eqref{TWlattice}, we directly get that $S'>0$ on $\R$ for any traveling wave solution. As $S$ is bounded, we get that $S(-\infty)=s_\infty$ exists. Since $I(\pm\infty)=0$, we infer that $I'(\pm\infty)=0$ since $c>0$ and as $I$ is a smooth profile. From the first equation we get that
\bqs
I(x)=c \frac{S'(x)}{\tau S(x)}.
\eqs
Injecting this expression into the second expression and integrating on the real line, we obtain
\begin{align*}
0 &= c\int_{-\infty}^{+\infty} S'(x) \left( 1- \frac{\eta}{\tau S(x)}\right)\md x + \lambda \int_{-\infty}^{+\infty}I(x-1)-2I(x)+I(x+1)\md x \\
&= c \int_{-\infty}^{+\infty} \Psi(S(x))' \md x = c \left( \Psi(s_0)- \Psi(s_\infty)\right),
\end{align*}
where $\Psi(v):=v-\frac{\eta}{\tau}\ln(v)$. As a consequence, $s_\infty$ is the unique positive real such that $\Psi(s_\infty)=\Psi(s_0)$ with $0<s_\infty<s_0$.

%
%

If we denote $\mathcal{I}(x):= \frac{1}{c}\int_x^{+\infty} I(z)\md z$, we obtain that
\bqq
-c \mathcal{I}'(x)=  s_0 \left( 1-e^{-\tau\mathcal{I}(x)}\right)-\eta \mathcal{I}(x) +\lambda (\mathcal{I}(x-1)-2\mathcal{I}(x)+\mathcal{I}(x+1)), \quad x\in\R,
\label{TWkpplikeD}
\eqq
with asymptotic conditions
\bqq
\mathcal{I}(-\infty)=\mathcal{I}_*, \quad \text{ and } \quad \mathcal{I}(+\infty)=0,
\label{TWkpplikeDlimit}
\eqq
with $\mathcal{I}'<0$ on $\R$. We remark that $\mathcal{I}_*= \frac{s_0-s_\infty}{\eta}=\frac{1}{\tau}\ln \frac{s_0}{s_\infty}$. The existence and uniqueness of monotone traveling wave solutions of \eqref{TWkpplikeD}-\eqref{TWkpplikeDlimit} is well-known. The existence was first proved in \cite{ZHH91} and then extended in \cite{CG02}. Uniqueness together with some refined asymptotic properties on the profile of the traveling waves can be found in \cite{CFG06}. We summarize these results in the following proposition from which our Theorem~\ref{thm4} easily follows.

\begin{prop}
Let us assume that $\mathscr{R}_0>1$. Then, there exists a unique (up to translation) monotone traveling wave profile solution of \eqref{TWkpplikeD}-\eqref{TWkpplikeDlimit} if and only if $c\geq c_*$.
\end{prop}

\section{Spreading properties on $\mathbb{T}_k$ with $k\geq2$}\label{secT}

In this section, we turn our attention to the spreading properties on $\mathbb{T}_k$ with $k\geq2$ for system \eqref{KPPtree}. Most of the proofs remain unchanged compared to the previous case on the lattice. Here, we will mostly highlight the key differences.

The proof of Theorem~\ref{thm5} follows the same strategy as for the proof of Theorem~\ref{thm1} by constructing constant sub and supersolutions to get the existence of a stationary solution. Note that the comparison principle applies to system \eqref{KPPtree} by direct application of the formalism developed in \cite{chen97}. Positivity of the stationary solution can also be established by contradiction. Assume that for some $n_0\geq 1$ we have $\mathcal{I}_{n_0}^\infty=0$. If $n_0=1$, the second equation of \eqref{KPPtree} gives
\bqs
0 = i_0+ \lambda (k+1) \mathcal{I}_{2}^\infty,
\eqs
and we have reached a contradiction since $i_0>0$. If now $n_0>1$, using the first equation of \eqref{KPPtree} gives necessarily that $\mathcal{I}_{n_0-1}^\infty=\mathcal{I}_{n_0+1}^\infty=0$. By induction, we must also have $\mathcal{I}_{1}^\infty=0$ which is impossible. 

Regarding the asymptotic behavior of $(\mathcal{I}_n^\infty)_{n\geq1}$ as $n\rightarrow+\infty$ when $\mathscr{R}_0\geq 1$, it is exactly the same as in the case on the lattice by noticing that $f<0$ on $(0,+\infty)$ in that regime. Thus, we  only have to treat the case $\mathscr{R}_0> 1$ and we distinguish two cases.

\begin{itemize}

\item Case $\lambda>\lambda_c$. As explained in the introduction, in this regime, one can find $\gamma>0$ such that
\bqs
\mathcal{D}(\gamma)=\eta(\mathscr{R}_0-1)+\lambda\left( e^{\gamma}-(k+1)+ke^{-\gamma}\right)<0.
\eqs
As a consequence, the sequence $(\overline{\mathcal{I}}_n^\infty)_{n\geq1}=(Ce^{-\gamma n})_{n\geq1}$ satisfies
\bqs
f( \overline{\mathcal{I}}_n^\infty)+\lambda\left(\overline{\mathcal{I}}_{n-1}^\infty-(k+1)\overline{\mathcal{I}}_n^\infty+k\overline{\mathcal{I}}_{n+1}^\infty\right)
\eqs
for each $n\geq2$ and $C>0$. We now select $C>0$ large enough such that
\bqs
f(Ce^{-\gamma})+i_0+\lambda C(k+1)(-1+e^{-\gamma})e^{-\gamma}<0,
\eqs
which is always possible since $f<0$ on $(\mathcal{I}_*,+\infty)$. Thus, $(\overline{\mathcal{I}}_n^\infty)_{n\geq1}=(Ce^{-\gamma n})_{n\geq1}$ is a supersolution and for all $n\geq 1$ we have $\mathcal{I}_n^\infty \leq C e^{-\gamma n}$ which implies that
\bqs
\mathcal{I}_n^\infty \longrightarrow 0 \text{ as } n \rightarrow+\infty.
\eqs
\item Case $0<\lambda<\lambda_c$. We introduce the sequence of shifts $\mathcal{I}_{n+\tau_m}^\infty$ with $\tau_m\in\N$ for $m\in \N$ and $\tau_m\rightarrow +\infty$ as $m\rightarrow+\infty$. Up to subsequences, it converges towards a nonnegative bounded sequence $(\widetilde{\mathcal{I}}_n)_{n\geq1}$ which satisfies 
\bqs
\left\{
\begin{split}
0&=f(\widetilde{\mathcal{I}}_n)+\lambda\left(\widetilde{\mathcal{I}}_{n-1}-(k+1)\widetilde{\mathcal{I}}_n+k\widetilde{\mathcal{I}}_{n+1}\right), \quad n \geq 2,\\
0&=f(\widetilde{\mathcal{I}}_1)+\lambda(k+1)\left(-\widetilde{\mathcal{I}}_1+\widetilde{\mathcal{I}}_2\right).
\end{split}
\right. 
\eqs
Note that by construction the stationary solution $(\mathcal{I}_n^\infty)_{n\geq1}$ is  a supersolution to the Fisher-KPP equation set on the homogeneous tree
\bqs
\left\{
\begin{split}
\mathcal{I}'_n(t)&=f( \mathcal{I}_n(t))+\lambda\left(\mathcal{I}_{n-1}(t)-(k+1)\mathcal{I}_n(t)+k\mathcal{I}_{n+1}(t)\right), \quad n \geq 2,\\
\mathcal{I}'_1(t)&=f(\mathcal{I}_1(t))+\lambda(k+1)\left(-\mathcal{I}_1(t)+\mathcal{I}_2(t)\right).
\end{split}
\right. 
\eqs
Using \cite{HH19}, we have that for $0<\lambda<\lambda_c$ any solution from a positive bounded  compactly supported initial converges locally uniformly as $t\rightarrow+\infty$ to the positive zero $f$ which is $\mathcal{I}_*>0$. By comparison, we have $\mathcal{I}_n^\infty\geq\mathcal{I}_*$ for each $n\geq1$. This shows that necessarily $\widetilde{\mathcal{I}}_n\geq\mathcal{I}_*$ and thus $\widetilde{\mathcal{I}}_n=\mathcal{I}_*$ since $f<0$ on $(\mathcal{I}_*,+\infty)$.
\end{itemize}

It remains to prove the uniqueness of the stationary solution. Once again, we distinguish between two cases.

\begin{itemize}
\item Case $\mathscr{R}_0>1$ and $0<\lambda<\lambda_c$. Let $\left(\mathcal{I}_n^\infty\right)_{n\geq1}$ and $\left(\mathcal{J}_n^\infty\right)_{n\geq1}$ be two positive, bounded stationary solutions to \eqref{KPPtree}. Since $\mathcal{I}_n^\infty,\mathcal{J}_n^\infty \geq \mathcal{I}_*$ and both $\mathcal{I}_n^\infty$ and $\mathcal{J}_n^\infty$ are bounded for each $j$, we have that
\bqs
\theta:=\underset{n\geq1}{\sup} ~ \frac{\mathcal{I}_n^\infty}{\mathcal{J}_n^\infty}>0,
\eqs
is a well-defined quantity. Assume by contradiction that $\theta>1$. Since $\underset{n\rightarrow+\infty}{\lim}\mathcal{I}_n^\infty=\underset{n\rightarrow+\infty}{\lim}\mathcal{J}_n^\infty=\mathcal{I}_*$, we have that the above supremum is a maximum attained at some $n_0\geq1$. Assume that first that $n_0>1$. We obtain by subtracting the equations evaluated at $n=n_0$ that  
\begin{align*}
0=f(\mathcal{I}_{n_0}^\infty)-\theta f(\mathcal{J}_{n_0}^\infty) &+\underbrace{\lambda\left( \mathcal{I}_{{n_0}-1}^\infty-2\mathcal{I}_{n_0}^\infty+\mathcal{I}_{{n_0}+1}^\infty\right)-\theta \lambda\left( \mathcal{J}_{{n_0}-1}^\infty-2\mathcal{J}_{n_0}^\infty+\mathcal{J}_{{n_0}+1}^\infty\right)}_{\leq 0}\\
&+\underbrace{\lambda(k-1)\left( \mathcal{I}_{{n_0}+1}^\infty-\mathcal{I}_{n_0}^\infty\right)-\theta \lambda (k-1)\left( \mathcal{J}_{{n_0}+1}^\infty-\mathcal{J}_{n_0}^\infty\right)}_{\leq 0}.
\end{align*}
This yields that
\bqs
\theta f(\mathcal{J}_{n_0}^\infty) \leq f(\mathcal{I}_{n_0}^\infty) = f(\theta \mathcal{J}_{n_0}^\infty),
\eqs
which is impossible by concavity of $f$. If now $n_0=1$, we obtain that
\bqs
0=f(\mathcal{I}_{1}^\infty)-\theta f(\mathcal{J}_{2}^\infty)+\underbrace{i_0(1-\theta)}_{<0}+\underbrace{\lambda(k+1)\left( -\mathcal{I}_{1}^\infty+\mathcal{I}_{2}^\infty\right)-\theta \lambda(k+1)\left( -\mathcal{J}_{1}^\infty+\mathcal{J}_{2}^\infty\right)}_{\leq 0},
\eqs
which gives $\theta f(\mathcal{J}_{1}^\infty)\leq f(\theta\mathcal{J}_{1}^\infty)$ which is impossible by concavity of $f$. As a consequence, we have proved that $\mathcal{I}_n^\infty\leq \mathcal{J}_n^\infty$ for each $n\geq1$. By reversing the role of the two solutions, we reach the conclusion.
\item Case $\mathscr{R}_0>1$ and $\lambda>\lambda_c$ or $\mathscr{R}_0\geq1$. In these two cases, any positive, bounded stationary solutions $\left(\mathcal{I}_n^\infty\right)_{n\geq1}$ and $\left(\mathcal{J}_n^\infty\right)_{n\geq1}$ to \eqref{KPPtree} asymptotically converge to zero as $n\rightarrow+\infty$. Take $\epsilon>0$ and define the sequence with $\mathcal{J}_n^{\infty,\epsilon}=\mathcal{J}_n^\infty+\epsilon$ which satisfies
\bqs
\left\{
\begin{split}
f( \mathcal{J}_n^{\infty,\epsilon})+\lambda\left(\mathcal{J}_{n-1}^{\infty,\epsilon}-(k+1)\mathcal{J}_n^{\infty,\epsilon}+k\mathcal{J}_{n+1}^{\infty,\epsilon}\right)<0, \quad n \geq 2,\\
f(\mathcal{J}_1^{\infty,\epsilon})+i_0+\lambda(k+1)\left(-\mathcal{J}_1^{\infty,\epsilon}+\mathcal{J}_2^{\infty,\epsilon}\right)<0.
\end{split}
\right. 
\eqs
Assuming by contradiction that $\mathcal{I}_n^\infty>\mathcal{J}_n^{\infty,\epsilon}$ somewhere and repeating the same argument as before, we end up with an inequality of the form $\theta f(\mathcal{J}_{n_0}^{\infty,\epsilon})\leq f(\theta \mathcal{J}_{n_0}^{\infty,\epsilon})$ for some $\theta>1$ and $n_0\geq1$. This is impossible since $f$ is concave and vanishes at the origin. As $\epsilon>0$ is arbitrary, we conclude that $\mathcal{I}_n^\infty\leq\mathcal{J}_n^{\infty}$ for all $n\geq1$ which ends the proof.
\end{itemize}

The proof of Theorem~\ref{thm7} is identical to as the proof of Theorem~\ref{thm2}, and we let to the reader.

Finally, we turn to the proof of Theorem~\ref{thm8} where we recall that we assume that $\mathscr{R}_0>1$ and $0<\lambda<\lambda_c$. We first construct a supersolution. We recall the definition of $c_*^k>0$ as 
\bqs
c_*^k:=\underset{\gamma>0}{\min}~  \frac{\eta \left( \mathscr{R}_0-1\right)+\lambda \left( e^{\gamma}-(k+1)+ke^{-\gamma} \right)}{\gamma},
\eqs
where the minimum is achieved at a unique value $\gamma_*>0$.  We introduce the sequence $\mathcal{H}_n(t):=\alpha e^{-\gamma_* (n-c_*t)}$ for $n\geq1$ for some constant $\alpha>0$ which is fixed large enough such that $\mathcal{I}_1(t)<\mathcal{H}_1(t)$  for all $t>0$. For $n\geq2$, we compute
\begin{align*}
\mathcal{H}_n'(t)-f(\mathcal{H}_n(t))-&\lambda\left( \mathcal{H}_{n-1}(t)-(k+1)\mathcal{H}_n(t)+k\mathcal{H}_{n+1}(t)\right)\\
&>\left(-c_*\gamma_* -f'(0)-\lambda \left( e^{\gamma_*}-(k+1)+ke^{-\gamma_*} \right) \right) \mathcal{H}_n(t) =0,
\end{align*}
thanks to the concavity of $f$. Thus $(\mathcal{H}_n(t))_{n\geq1}$ is a supersolution and we can deduce that $\mathcal{I}_n(t)<\mathcal{H}_n(t)$ for all $t>0$ and $n\geq1$. As a consequence, for each $c>c_*^k$, we have
\bqs
\underset{t\rightarrow+\infty}{\lim} \left( \underset{n \geq ct}{\sup}\left| \mathcal{I}_n(t) \right| \right) =0.
\eqs

To conclude the proof of the theorem we let $c\in(0,c_*^k)$ and a consider a sequence $(t_m)_{m\in\N}$ such that $t_m\rightarrow+\infty$ as $m\rightarrow+\infty$ and a sequence $(n_m)_{m\in\N}$ in $\N$ such that $|n_m|\leq (c_*^k-c)t_m$. If $(n_m)_{m\in\N}$ is bounded, we know from Theorem~\ref{thm7} that $\mathcal{I}_{n_m}(t_m)-\mathcal{I}_{n_m}^\infty\rightarrow 0\rightarrow 0$ as $m\rightarrow+\infty$ by local uniform convergence. Suppose that up to subsequences $(n_m)_{m\in\N}$ diverges. We recall that the solution $\left(\mathcal{I}_n(t)\right)_{n\geq1}$ is a supersolution of the Fisher-KPP equation set on the homogeneous tree for which spreading occurs with the asymptotic speed $c_*^k$ thanks to the results of \cite{HH19}. We infer that
\bqs
\underset{m\rightarrow+\infty}{\liminf} \left( \mathcal{I}_{n_m}(t_m)-\mathcal{I}_{n_m}^\infty \right) \geq \mathcal{I}_*-\mathcal{I}_*=0.
\eqs
Next, we let $\gamma_0>0$ be the unique positive zero of  $\varphi_k(\gamma):=f'(\mathcal{I}_*)+\lambda \left( e^{\gamma}-(k+1)+ke^{-\gamma}\right)$. We define the sequence $\mathcal{H}_n:=\mathcal{I}_*+ \beta e^{-\gamma_0 n}$ for $n\geq1$. We readily remark by concavity of $f$ that for each $n\geq2$, we have
\bqs
f(\mathcal{H}_n)+\lambda \left( \mathcal{H}_{j-1}-(k+1)\mathcal{H}_n+k\mathcal{H}_{n+1}\right)<\left( f'(\mathcal{I}_*)+\lambda \left( e^{\gamma_0}-(k+1)+e^{-\gamma_0}\right) \right)\beta e^{-\gamma_0 n}=0.
\eqs
 Then, we select $\beta>0$ large enough, so that $\mathcal{I}_1(t)<\mathcal{H}_1$ for all $t>0$, which is possible since $\mathcal{I}_1(t)$ is bounded. Hence, we obtain by comparison that $\mathcal{I}_n(t)<\mathcal{H}_n$ for each $n\geq1$ and $t>0$. It follows that
\bqs
\underset{m\rightarrow+\infty}{\limsup}~ \left( \mathcal{I}_{n_m}(t_m) -\mathcal{I}_{n_m}^\infty\right)\leq \underset{m\rightarrow+\infty}{\limsup}~ \left(\mathcal{I}_*-\mathcal{I}_{n_m}^\infty+\beta e^{-\gamma_0 n_m} \right) =0.
\eqs
This concludes the proof of the theorem.

\section{Discussion}\label{disc}

In this paper, we investigated the spreading properties of a SIR model set on homogeneous trees of degree $k$ and distinguished between the case $k=1$ which reduces to the usual lattice of integers $\Z$ and the case $k\geq2$. By means of a classical transformation \cite{A77}, the cumulated density of infected individuals is shown to satisfy a non-homogeneous Fisher-KPP type equation which allows us to extend previous results known in the continuous case \cite{DG14,BRR21} to our discrete setting. We could derive some rather precise properties of our model which we summarize here.

First, we proved that the non-homogeneous Fisher-KPP equation satisfied by the cumulated density of infected individuals admits a unique bounded positive stationary solution which is a global attractor for the dynamics of this non-homogeneous Fisher-KPP equation. When the graph is the lattice $\Z$, the asymptotic behavior of the stationary solution depends on the basic reproduction $\mathscr{R}_0$, and there is a dichotomy. When $\mathscr{R}_0\leq1$, this steady state tends to zero at infinity, whereas if $\mathscr{R}_0>1$ it asymptotically converges towards a positive constant. In that case, we have propagation of the epidemic in the lattice and we further quantified the asymptotic speed of spreading as a function of the parameters of the model. In the regime of strong diffusion within the lattice, we recover the asymptotic speed of propagation in the continuum case. We have further characterized this speed as the threshold for the existence of traveling wave solutions to our system. When the graph is a homogeneous tree of degree $k\geq2$, we have also demonstrated  that the non-homogeneous Fisher-KPP equation satisfied by the cumulated density of infected individuals admits a unique bounded positive stationary solution which is a global attractor for the dynamics of this non-homogeneous Fisher-KPP equation. When $\mathscr{R}_0\leq1$, this steady state tends to $0$ at infinity.  Now when $\mathscr{R}_0>1$, we have proved the existence of a threshold on the strength of interactions, which depends on the degree of the tree and the basic reproduction number. Below this threshold, the steady state asymptotically converges to a positive constant and an epidemic can spread in the tree, and above the threshold it asymptotically converges to zero. Similarly to the case on the lattice, we managed to quantify the asymptotic speed of spreading as a function of the parameters of the model in the regime where spreading is possible. This asymptotic speed of propagation reaches a maximum at a critical value of the strength of interactions which can be explicitly computed. The maximal spreading speed, defined as the speed at this critical value, scales as $1/\ln k$ indicating that the higher the degree of the tree the smaller the asymptotic spreading spreed.

Our model sheds light on the effect of networks structure on the propagation of epidemics in the specific case of homogeneous trees. It allows us to explain some counterintuitive phenomena. Thanks to its relative simple structure, we could derive closed form formula and carry a somehow complete mathematical analysis. We see this work as a first step towards a more systematic understanding of spreading phenomena in more realistic and practical networks, such as transportation networks for example \cite{BH13}. In a very recent work \cite{BF21}, we proposed a new model that describes the dynamics of epidemic spreading on connected graphs. Our model consists in a PDE/ODE system where at each vertex of the graph we have a standard SIR model and connections between vertices are given by heat equations on the edges supplemented with Robin like boundary conditions at the vertices modeling exchanges between incident edges and the associated vertex. Under some appropriate scaling assumptions, the model of the present paper can be seen as the limit of the PDE/ODE model from \cite{BF21}. One of our objective will be to understand how the spreading properties analyzed here can be transposed to this extended PDE/ODE model.

\section*{Acknowledgements} 

This works was partially supported by Labex CIMI under grant agreement ANR-11-LABX-0040. G.F. acknowledges support from an ANITI (Artificial and Natural Intelligence Toulouse Institute) Research Chair.

\end{document}